# Two-stage battery recharge scheduling and vehicle-charger assignment policy for dynamic electric dial-a-ride services

**Tai-Yu Ma**

Department of Urban Development and Mobility, Luxembourg Institute of Socio-Economic Research (LISER), Esch-sur-Alzette, Luxembourg

# Abstract

Coordinating the charging scheduling of electric vehicles for dynamic dial-a-ride services is challenging considering charging queuing delays and stochastic customer demand. We propose a new two-stage solution approach to handle dynamic vehicle charging scheduling to minimize the costs of daily charging operations of the fleet. The approach comprises two components: daily vehicle charging scheduling and online vehicle–charger assignment. A new battery replenishment model is proposed to obtain the vehicle charging schedules by minimizing the costs of vehicle daily charging operations while satisfying vehicle driving needs to serve customers. In the second stage, an online vehicle–charger assignment model is developed to minimize the total vehicle idle time for charges by considering queuing delays at the level of chargers. An efficient Lagrangian relaxation algorithm is proposed to solve the large-scale vehicle-charger assignment problem with small optimality gaps. The approach is applied to a realistic dynamic dial-a-ride service case study in Luxembourg and compared with the nearest charging station charging policy and first-come-first-served minimum charging delay policy under different charging infrastructure scenarios. Our computational results show that the approach can achieve significant savings for the operator in terms of charging waiting times (–74.9%), charging times (–38.6%), and charged energy costs (–27.4%). A sensitivity analysis is conducted to evaluate the impact of the different model parameters, showing the scalability and robustness of the approach in a stochastic environment.

# Introduction

Electric vehicle technology has gained increasing interest amongst policymakers, the general public, and the automotive industry in response to worldwide directives to reduce $CO_2$ emissions. Transport network companies (TNCs) such as Lyft and Uber have begun deploying battery electric vehicles (called EVs hereafter) in their fleet to reduce operating costs and promote green mobility [1]. Research on the electrification of ride-hailing services in the USA has shown that TNCs need to recharge e-fleets several times a day and rely primarily on DC fast chargers to minimize charging times [2]. As charging EVs with high-power charging (>22kW) is much more expensive than residential electricity prices, such operations may significantly increase the operator's charging cost by about 25% [3]. Furthermore, due to the higher installation cost of high-power charging, most public chargers are limited to Level 2 chargers [4]. With limited battery range, a vehicle's within-day charging becomes a primary challenge given the limited number of charging stations and relatively long charging time. For example, an 80% charge needs about 50 minutes using a 50 kW DC fast charger for a Volkswagen Golf with a 300 kilometer range [4]. Additionally, with the increased number of electric vehicles in the fleet and the relatively limited number of public and private charging spots, the likelihood that accessible charging stations will be temporarily unavailable will soon become an issue. Uncoordinated charging operations might result in higher queueing delays, reducing the availability of vehicles to serve customers, and an increase in total system operating costs. However, existing studies mainly focus on static EV routing problems under charging infrastructure constraints, whereas research on online charging scheduling under stochastic demand is still limited [5]. For this purpose, we propose an online charging scheduling model for dynamic dial-a-ride services to minimize the total charging delays and costs of the fleet under charging infrastructure constraints and stochastic customer demand.

The challenge of charging scheduling for electric dynamic dial-a-ride services under uncertainty involves several dimensions. First, under stochastic customer arrivals, vehicle driving patterns are stochastic, which impacts vehicle charging demand in space and time. Second, given limited charging facility resources and the stochastic charging demand of other EVs, there might be queuing delays at charging stations. How to efficiently coordinate the charging demand of the e-fleet while considering vehicle's driving needs and charging station capacity constraints? Third, given heterogeneous charging powers and space-time differentiated charging prices, how should operators decide when, how much, and where to charge vehicles such that the overall charging costs and queuing delays are minimized? To address these challenges, we propose an online charging scheduling model by anticipating the future energy needs of vehicles to minimize the total charging delays and charging costs of the fleet of EVs for dynamic dial-a-ride services.

The remainder of the paper is organized as follows. Section 2 reviews related literature on EV charging strategies in a stochastic environment. Section 3 proposes a two-stage EV recharging policy for dynamic dial-a-ride services. We first derive an optimal charging plan for each individual vehicle based on its historical driving patterns, the price of electricity, and expected queuing delays at charging stations. Then an optimal charging station assignment model is proposed to minimize total vehicle charging times and queuing delays based on the current system state. An efficient solution heuristic based on the Lagrangian relaxation (LR) method is proposed for dealing with large-scale instances and allowing real-time operations. In Section 4, we conduct a realistic case study for a dynamic dial-a-ride service in Luxembourg to evaluate the performance of the proposed methodology in a stochastic environment. The impact of different model parameters on system performance is analyzed. Finally, conclusions are drawn and future extensions are discussed.

# Related work

Coordinating EV charging demand in order to reduce its impact on the electric grid has been studied in recent years [6, 7]. However, these charging scheduling models are mainly from a private EV owner perspective, which mainly involves recharging EVs at home or at the workplace once or twice a day. For a fleet operator, the charging optimization strategy is different from that of private EVs in both time (charging operations during the day) and scale (the vehicle fleet needs to be recharged several times per day). Charging coordination strategy needs to consider various factors of uncertainty: stochastic vehicle driving patterns, queuing delays at charging stations, charging price variations, charging infrastructure capacity constraints, and customer inconvenience due to recharging EVs. Hu et al. [8] classify three EV charging control strategies for fleet operators: centralized control, transactive control, and price control. Centralized control assumes the operator directly schedules EV recharging operations via real-time communication. Transactive control is a kind of distributed control mechanism to achieve supply–demand equilibrium in an electricity market. Price control relies on a dynamic electricity pricing design to regulate electricity supply and demand disequilibrium for EVs. The strategy of coordinating spatio-temporal supply-demand mismatch to enhance the system efficiency has also been studied for the collaborative logistics problems [9, 10].

Iacobucci [11] pointed out that studies for the design of charging strategies for shared EVs are still limited. The authors propose a two-layer model predictive control strategy for relocating and charging shared autonomous electric vehicles (SAEVs). Several charging optimization and idle vehicle relocation models have been proposed for electric car-sharing systems [12-15]. For dynamic dial-a-ride service using EVs, a number of works have proposed mathematical models for optimizing shared electric autonomous vehicles operations. For example, Zhang and Chen [16] proposed a charging optimization strategy to balance the charging demand of SAEVs in a high-priced electricity period to reduce total charging cost. The battery levels of individual vehicles are first sorted, and then individual vehicles with low battery levels are set to recharge. The number of concurrent charging SAEVs is regulated by the ratio between the energy demand of SAEVs, the available number of chargers for SAEVs, and SAEV recharge rates. Queuing delays and charging station assignment are not explicitly considered. Bongiovanni [17] proposed a two-phase metaheuristic to solve the dynamic dial-a-ride problems using electric autonomous vehicles. The proposed approach first solves a static dial-a-ride problem under battery constraint to initiate vehicles' routing plans. New requests are then inserted into vehicles' planned routes to minimize a weighted objective function by considering both operational costs, customer inconvenience, and demand. A two-phase metaheuristic is proposed to find good EV vehicle routing solutions. Zalesak and Samaranayake [18] proposed a batch-optimization framework based on the shareability network concept for ride-pooling using EVs. New requests are assigned first to vehicles under current charging schedule constraints. Upon new assigned requests, a charging planning model is used to update the charging schedules of vehicles. A more realistic battery charging model considering the nonlinear relationship of battery state of charge and charging time is used. Some recent studies formulated the dynamic ridesharing problems using EVs as a Markov decision process and proposed approximate dynamic programming approaches to maximize the profit of operators [19, 20].

Several studies propose mathematical models to evaluate the impact of deploying e-taxis or SAEVs on the level of service. The charging strategy is mainly based on a full-charge policy to recharge EVs to a maximum level at nearest charging stations whenever an EV's battery level is lower than a threshold [21, 22]. Tian et al. [23] proposed a real-time charging station recommendation system for e-taxis based on the historical driving patterns of vehicles. When receiving a vehicle's charging request, the recommendation system suggests the charging station with the least total access time and waiting time when arriving at the charging station, according to the order of received requests. The results show that the proposed recommendation system could significantly reduce vehicle's waiting times compared to the nearest charging station assignment policy. From an individual taxi driver perspective, this first-come-first-served policy is most beneficial for each new charging request. However, from a fleet charging management perspective, charging operations can be further optimized by coordinating the charging demand over non-rush hours and allowing partial recharging to reduce taxis' idle time. For example, Yuan et al. [24] proposed an e-taxi charging scheduling model under a receding horizon control framework allowing partial recharging to minimize taxi fleet idle time under dynamic taxi demand. The results suggest that partial charging allows for reducing vehicle waiting times and increasing the number of available taxis in rush hour. However, the considered charging infrastructure is assumed to be homogeneous and charging station assignment is not optimized to minimize total queuing delays.

The partial recharge policy raises the issue of what battery levels are necessary to satisfy a vehicle's driving needs and how to determine the optimal charging plans for vehicles based on individual vehicle's driving patterns. Iversen et al. [25] proposed a model to optimize the charging level plan of an individual PEV based on individual vehicles' historical driving patterns. The problem is considered as a stochastic dynamic programming problem to minimize the total charging cost while satisfying the vehicle's energy needs for driving. An inhomogeneous Markov model is fitted by using individual vehicles' stochastic driving patterns to estimate the state transition probability from being idled to a driving state. A summary of existing studies on charging policies for on-demand shared mobility services is shown in Table 1.

Dynamic dial-a-ride problems using e-fleets present a more complex environment for managing charging operations under stochastic customer demand and charging capacity constraints. To the best of our knowledge, the current state of the art has not fully addressed these issues to minimize total e-fleet idle times and charging cost under dynamic customer demand and queuing delays at charging stations.

The main contributions of the current work are summarized as follows.

1) We propose a two-stage approach to handle the vehicle charging scheduling problem for dynamic dial-a-ride services using EVs to minimize the daily charging operational costs and delays of the fleet. A first vehicle charging scheduling model is formulated as a battery recharge problem under uncertainty to minimize the total charging operational costs by considering vehicle probabilistic driving needs, expected charging delays, and charging costs.

2) A new online vehicle–charger assignment model is proposed as a mixed-integer optimization problem to minimize the total vehicle idle times for recharges considering queuing delays at the level of the chargers. A Lagrangian relaxation algorithm is developed and tested on large-scale test instances. The computational results show that the LR algorithm can obtain near-optimal solutions within a couple of seconds/minutes for median-/large-sized problems.

3) A realistic dynamic dial-a-ride service case study is implemented to assess the performance of the proposed solution. The results show that significant savings in terms of charging delay, charging time, and cost can be achieved compared to the state-of-the-art nearest charging station policy and minimum charging delay policy.

**Table 1. Summary of charging policies for on-demand shared mobility services.**

| Studies | System | Charging policy features |
|---|---|---|
| Bischoff and Maciejewski [21]; Chen et al. [22] | e-taxi | On-need policy to assign a vehicle to the nearest charging station for recharge whenever the vehicle battery level is lower than a threshold. |
| Iacobucci [11] | SAEV | Consider the dynamic electricity price for scheduling vehicle charges in smart grids. A two-layer model predictive control approach is proposed to optimize vehicle charging scheduling over a longer timeframe. Congestion at charging stations is not considered. |
| Tian et al. [23] | e-taxi | Consider the inference of electric taxi states based on historical taxi charging patterns and position tracking. Uses the first-come-first-served policy for charging station allocation whenever vehicle charging intention is identified. |
| Yuan et al. [24] | e-taxi | Propose a zone-based charging station allocation policy to minimize vehicle idle times for recharge. Partial recharge is allowed without queuing delay consideration. |
| Ma et al. [12]; Pantelidis et al. [13] | carsharing | Static carsharing vehicle charge scheduling and relocation based on the facility location model. Stochastic demand and queuing delays are considered to meet customer demand. |
| Roni et al. [15] | carsharing | A capacitated facility location model is proposed for optimal charging station allocation on a time-space network to minimize total travel and waiting times of charging operations. |
| Folkestad [14] | carsharing | Propose a static carsharing vehicle charging scheduling and repositioning model to satisfy charging needs with minimal vehicle relocation costs. |
| Zhang and Chen [16] | SAEV | Propose a probabilistic rule for charging station allocation to regulate charging demand (number of vehicles sent to charge) and supply (number of available chargers). Considers electricity price variation to minimize charging costs without queuing delay considerations. |
| Rinaldi et al., [26]; Wang et al. [27] | e-bus | Propose a static electric bus charging and route planning model to minimize the total operational costs of the fleet. |
| Mkahl et al. [28] | fleet of electric vehicles | Propose a linear programming model for charging station allocation to keep a vehicle's battery at its highest possible level when arriving at a charging station. Full-charge policy without queuing delay considerations. |
| Lu et al. [29] | e-taxi | Propose a multi-commodity network flow model on a space-time network for a mixed fleet of EVs and gasoline vehicles. Travel requests are deterministic and known. No charge queuing delay consideration. |
| Al-Kanj et al. [19]; Yu et al. [20] | ridesharing | Modeling dynamic electric ridesharing problems as a Markov decision process and propose approximate dynamic programming approaches to maximize operator's profit. |

Remark: SAEV: shared autonomous electric vehicle.

# Methodology

We consider a dynamic dial-a-ride problem in which a TNC operates a fleet of homogeneous EVs to pick up and drop off customers. Ride requests arrive stochastically and are accepted/rejected on short notice. The fleet of EVs is assumed charged at a certain level (80% or more) at the beginning of day to ensure a good battery lifespan [16, 21]. A limited number of charging points are available in the service area to allow EVs to recharge. A dispatching center is equipped with a dedicated management platform with real-time information on vehicles (location and battery level) and charging station status (i.e., number and characteristics of chargers, and charging schedules of EVs at the location) [8]. The operator dispatches vehicles to pick up customers according to a designed vehicle routing and dispatching policy (described in Section 4.1). An EV's battery energy level is monitored in real-time and communicated to the dispatch center. Given stochastic vehicle driving patterns, uncertain charging demand from other EVs, and capacitated charging infrastructure, it is not possible to obtain the exact charging plans of vehicles (when, where, and how much energy to recharge each vehicle) in advance. The dynamic dial-a-ride charging scheduling problem is to design an online charging policy under these uncertainty factors to minimize total charging delays and costs of the e-fleet over the planning horizon (one day) under a stochastic environment.

For this purpose, we discretize the planning horizon into a set of charging decision epochs and decompose the decision process into two stages. In the first stage, we determine in advance the optimal vehicle charging schedules (when and how much energy to charge) for each epoch based on the historical driving patterns of vehicles and the expected time-dependent queueing delays at charging stations. The problem is formulated as a single-vehicle battery recharge problem to minimize total charging delays and costs while satisfying vehicle driving needs. At this stage, the specific charging station location assignment is not considered and waiting times to be served at charging stations are based on historical information. Given the charging schedules obtained in the first stage, the second stage determines the optimal vehicle–charger assignment by solving the charging station assignment problem to minimize total charging delays. Our computational study shows that the proposed methodology can effectively reduce total charging delays and system operational costs in a stochastic environment. The two-stage battery recharge scheduling framework is shown in Fig 1.

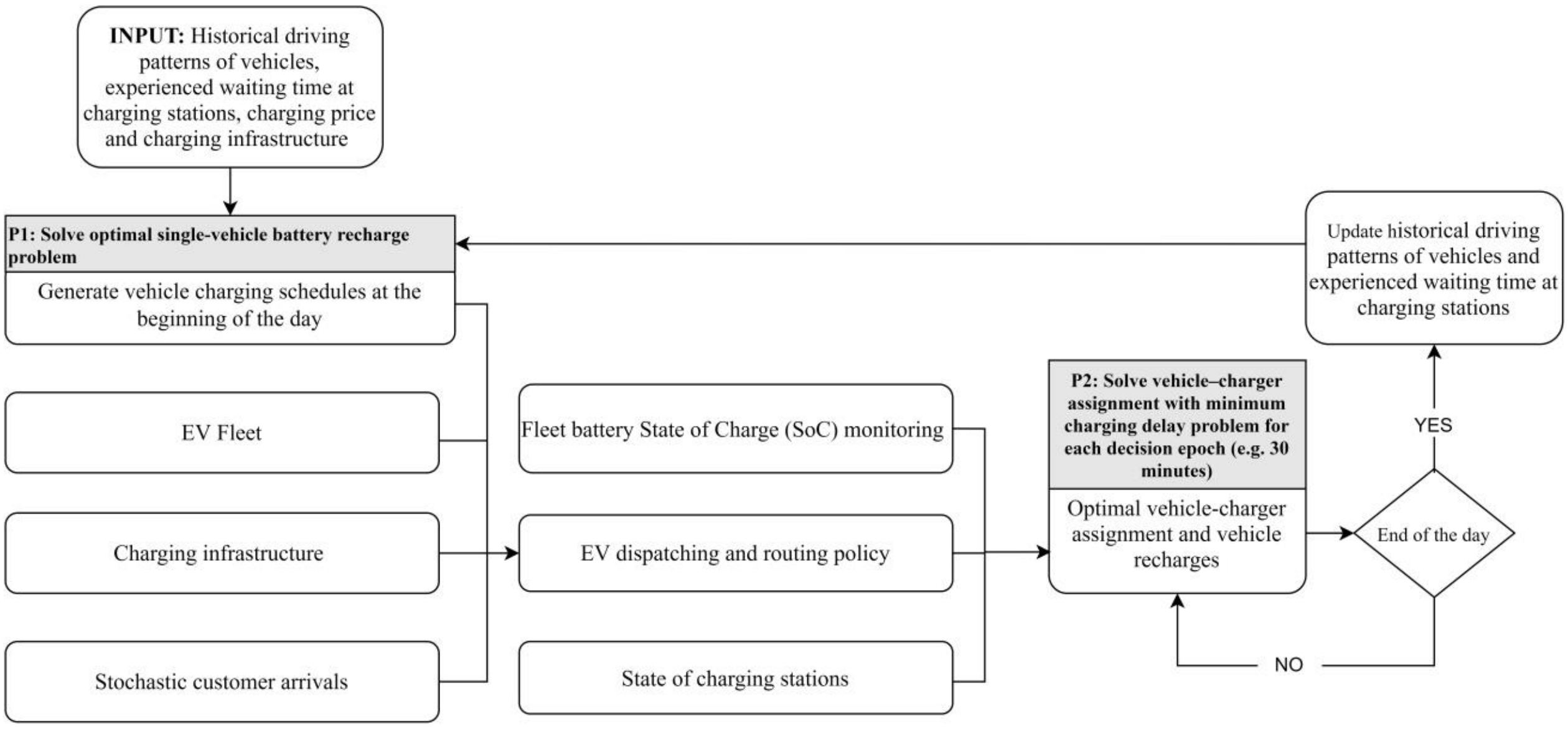

**Fig 1. Two-stage battery recharge scheduling framework.**

# Optimal vehicle charging schedules under stochastic driving patterns of vehicles

Notation

| | |
|---|---|
| $h$ | Index of charging decision epochs, $h \in H = \{1,2,\dots,\lvert\mathrm{H}\rvert\}$ |
| $e_h$ | Energy level (state of charge) of a vehicle at the beginning of decision epoch $h \in H$(kWh) |
| $d_h$ | Expected energy consumption in epoch $h$ based on vehicles' historical driving patterns (kWh) |
| $e_{max}$ | Battery capacity or the allowed maximum energy level of vehicle (kWh) |
| $e_{min}$ | Reserved energy level of a vehicle, $e_{min} = 0.1B$ (kWh) |
| $\vartheta$ | Energy price (euro/kWh) |
| $\varphi$ | Charging rate of chargers (kW/min.) |
| $\rho$ | Average gains per minute travelled (euro/min.) |
| $v$ | Vehicle speed (km/min.) |
| $\Delta$ | Time interval between any two consecutive epochs (min.) |
| **Decision variables** | |
| $u_h$ | Amount of charged energy in epoch $h$ (kWh) |
| $y_h$ | 1 if a vehicle is recharged in epoch $h$, and 0 otherwise. |

Given that the considered single-vehicle battery recharge problem has an intrinsic multi-stage decision-making nature in a stochastic environment, the problem is decomposed into a sequence of simpler one-stage decision problems. The entire planning period is discretized into a set of recharge decision epochs, $H = \{1,\dots,|H|\}$, with a time interval $\Delta$. The system state is the battery level $e_h$ of a vehicle at the beginning of each epoch $h$. The decision (control) variable is the amount of energy to charge $u_h$ for each epoch $h$. A cost function is associated with the charging decision, which depends on the cost of the energy charged and the opportunity cost of the unavailability of vehicles to serve customers. For simplification, a linear charging efficiency is assumed that charging time is equal to the amount of charge divided by the charging efficiency of chargers. Given the average energy consumption from historical driving patterns in each epoch, the system state is updated at the end of each epoch and an optimal recharge policy can be derived over the planning horizon. The optimal EV battery recharge problem is formulated as follows.

<u>P1: Optimal single-vehicle battery recharge problem</u>

$$\mathrm{Min} \sum_{h=1}^{|H|} [\vartheta u_h + (\bar{c} + \omega_h) y_h] \tag{1}$$

subject to

$$e_{h+1} = e_h + u_h - d_h, \quad \text{for } h = 1,\dots,|H| \tag{2}$$

$$e_h + u_h \geq d_h + e_{min}, \quad \text{for } h = 1,\dots,|H| \tag{3}$$

$$u_h \leq M y_h, \quad \text{for } h = 1,\dots,|H| \tag{4}$$

$$e_1 = e_{max} \tag{5}$$

$$e_{min} \leq e_h \leq e_{max}, \quad \text{for } h = 1,\dots,|H| \tag{6}$$

$$0 \le u_h \le u_{max}, \quad \text{for } h = 1, \dots, |H| \tag{7}$$

$$y_h \in \{0,1\}, \text{for } h = 1, \dots, |H| \tag{8}$$

The objective function (1) is to minimize the total cost of charging operations of the vehicle over the planning horizon. The first term relates to the cost of the amount of charged energy. The second term $\bar{c}$ is a fixed operational cost, which takes into account the set-up cost and the average access cost (energy consumption of travel distance to reach the charging stations). The last term $\omega_h$ is the opportunity cost of not being able to serve customers during recharging operations in decision epoch $h$, defined as Eq. (9):

$$\omega_h = \left(T_h^D + T_h^W\right)\rho \tag{9}$$

where $T_h^D$ is the expected duration of a vehicle being in use in $h$, estimated by the historical driving patterns of vehicles. $T_h^W$ is the expected waiting time to be served at charging stations in $h$. We estimate $T_h^D = \int_0^\Delta p_t^D \, dt$ and $T_h^W = \int_0^\Delta p_t^W \, dt$, where $p_t^D$ and $p_t^W$ are the probability of driving and waiting status, respectively, for a vehicle at time $t$. $\rho$ is the weight introduced to convert a vehicle's unavailable time due to charging operations to the gain loss based on a vehicle's average earnings. Eq. (2) is the state transition function describing the evolution of energy levels of the vehicle in each epoch. The expected energy consumption $d_h$ is defined as the total driving distance in epoch $h$ divided by the energy efficiency of vehicles; $d_h = v\Delta\mu$, where $\mu$ is the driving efficiency of vehicles (kWh/km). Eq. (3) indicates that the total energy level after recharge needs to be no less than the energy demand plus a minimum reserve energy $e_{min}$. Eq. (4) ensures that the amount of charged energy is non-negative when the vehicle goes to charge. Eq. (5) is the initial battery level of the vehicle. Eq. (6) states the upper and lower bounds of the energy level at the beginning of each epoch. Eq. (7) states that the amount of energy that can be charged for each epoch is bound by $u_{max} = \Delta\varphi$. Note that $u_{max}$ is an upper bound that the amount of energy can be charged on an epoch based on the fastest charger in a study area. The different charging rates of chargers are considered in the second stage for the vehicle-charger assignment to minimize the overall charging operational cost.

Problem P1 can be efficiently solved by the dynamic programming approach using the backward induction algorithm [30] or by a standard commercial mixed-integer optimization solver. Note that in a stochastic environment, each vehicle has different driving patterns during the planning horizon, so P1 needs to be solved for each vehicle to obtain the appropriate charging plans based on its historical driving patterns. The interactions between different vehicles are considered in the vehicle dispatching policy (described in Sect. 4.1). An illustrative example is given in Appendix A to illustrate the model property and the total charging operational cost savings compared with a reference on-need charging policy.

# Charging station assignment under charging capacity constraints

*Notation*

| | |
|---|---|
| $I$ | Set of vehicles to be recharged at the beginning of a recharging epoch $h$ (index $h$ is dropped) |
| $J$ | Set of chargers in a studied area |
| $t_{ij}$ | Travel time from the location of vehicle $i$ to that of charger $j$ |
| $d_{ij}$ | Travel distance from the location of vehicle $i$ to that of charger $j$ |
| $e_i$ | Energy level of vehicle $i$ at the beginning of epoch $h$ (index $h$ is dropped) |
| $e_i^*$ | Energy level of vehicle $i$ after recharge at the end of epoch $h$, determined by the charging plan from P1 (index $h$ is dropped). |
| $t_j^A$ | Time until which a charger $j$ is occupied by other vehicles from the beginning of epoch $h$ (index $h$ is dropped) |
| $\mu$ | Driving efficiency of vehicles (kWh/km) |
| $\varphi_j$ | Charging rate of charger $j$ (kW/min.) |
| $M$ | Large positive number |
| **Decision variable** | |
| $X_{ij}$ | Vehicle $i$ is assigned to charger $j$ for recharge if $X_{ij} = 1$, and 0 otherwise |
| $Y_{ij}$ | Amount of energy recharged at charger $j$ for vehicle $i$ |
| $W_{ij}$ | Artificial variable representing the waiting time of vehicle $i$ at charger $j$ |

For the second stage, the problem is to assign vehicles to chargers for each epoch $h \in H$ based on the charging schedules obtained beforehand. Given an epoch $h$, the charging delay of a vehicle at a charging station is defined as the sum of access time (travel time) to the charging station, the waiting time to be served at the charging station, and the total charging time. The problem is formulated as a mixed-integer optimization problem to minimize charging delays given that the capacitated charging infrastructure is a multi-server queuing system. Different from existing charging-station-based capacity constraints (the number of vehicles assigned to a charging station cannot exceed the number of chargers at that station [13, 15]), we consider each charger explicitly to account for the exact waiting time of a vehicle when arriving at a charger at time $t$ and the charging power of each individual charger.

The one-stage optimal charging station assignment model is formulated as follows. The problem is solved for each decision epoch $h \in H$.

P2: Vehicle–charger assignment with minimum charging delay problem

$$\min Z = \sum_{i \in I} \sum_{j \in J} t_{ij} X_{ij} + \theta_1 \sum_{i \in I} \sum_{j \in J} Y_{ij}/\varphi_j + \theta_2 \sum_{i \in I} \sum_{j \in J} W_{ij} \quad (10)$$

subject to

$$\sum_{j \in J} X_{ij} = 1, \qquad \forall i \in I \quad (11)$$

$$\sum_{i \in I} X_{ij} \leq 1, \qquad \forall j \in J \quad (12)$$

$$e_{min} \leq e_i - \mu d_{ij} X_{ij} + M(1 - X_{ij}), \forall i \in I, j \in J \quad (13)$$

$$e_i^* \leq Y_{ij} + e_i - \mu d_{ij} X_{ij} + M(1 - X_{ij}), \forall i \in I, j \in J \quad (14)$$

$$Y_{ij} \leq M X_{ij}, \forall i \in I, j \in J \quad (15)$$

$$t_j^A - t_{ij} X_{ij} - M(1 - X_{ij}) \leq W_{ij}, \forall i \in I, j \in J \quad (16)$$

$$X_{ij} \in \{0,1\}, \forall i \in I, j \in J \quad (17)$$

$$Y_{ij} \geq 0, \forall i \in I, j \in J \quad (18)$$

$$W_{ij} \geq 0, \forall i \in I, j \in J \quad (19)$$

The objective function minimizes the total weighted time of charging operations, including total travel time to arrival at charging stations, recharging time, and waiting time at each charging connector. $\theta_1$ and $\theta_2$ are the weights introduced to account for the trade-off between these elements. Constraints (11) and (12) ensure that each vehicle can be assigned to one charger and that each charger can be plugged in to at most one vehicle; constraint (13) guarantees that the remaining battery level of a vehicle when arriving at a charging station is no less than a pre-defined reserve level, e.g., 10%-20% of battery capacity. Constraint (14) states that the energy level after recharge must be no less than the planned level after recharge from P1. Constraint (15) ensures that the amount of recharged energy is non-negative when the vehicle is assigned to a charger for recharge. Constraint (16) calculates the waiting time to be served for vehicle $i$ when arriving at the location of charger $j$. Note that Eqs. (11) and (12) are suitable for the situation where the number of vehicles is no more than that of charges ($|I| \leq |J|$). In case of $|I| > |J|$, constraints (11) and (12) are replaced by (20) and (21), respectively.

$$\sum_{j \in J} X_{ij} \leq 1, \qquad \forall i \in I \tag{20}$$

$$\sum_{i \in I} X_{ij} = 1, \qquad \forall j \in J \tag{21}$$

We refer to the problem of Eqs. (10)–(19) as P2, and to Eqs. (10) and (13)–(21) as P2J. The above vehicle–charger assignment problem is a variant of the generalized assignment problem with additional constraints. We propose a heuristic based on the LR method to solve it for large instances in order to obtain efficiently near-optimal solutions for real-time applications.

# Proposed Lagrangian relaxation algorithm

The LR method is a widely-used methodology for solving mixed-integer optimization problems [31, 32]. This method first solves an LR problem by relaxing complicated constraints to obtain a lower-bound (LB) solution. As the LB solution is likely infeasible for the original problem, a problem-specific repair procedure needs to be developed to find a feasible solution, providing an upper bound (UB) to the original problem. Afterwards, the Lagrangian multiplier is updated to maximize the LB. The above steps are repeated until no improvement can be found or the maximum iteration is reached. For the P2 problem, it is not difficult to find that we can reformulate it by removing $Y$ and $W$ as follows:

$$\min \mathrm{Z} = \sum_{i \in I} \sum_{j \in J} (t_{ij} + \theta_2 \max(t_j^A - t_{ij}, 0)) X_{ij} + \theta_1 \sum_{i \in I} \sum_{j \in J} \frac{1}{\varphi_j} (e_i^* - e_i + \mu d_{ij}) X_{ij} \tag{22}$$

subject to (11)–(13) and (17).

We relax constraint (12) and use a non-negative Lagrangian multiplier $\lambda_j, \forall j \in J$ to penalize the non-satisfaction of this constraint in the objective function. The LR problem can be written as follows.

$$\min Z_{LB}(\boldsymbol{\lambda}) = \sum_{i \in I} \sum_{j \in J} \left[ t_{ij} + \lambda_j + \theta_2 \max(t_j^A - t_{ij}, 0) + \frac{\theta_1}{\varphi_j} (e_i^* - e_i + \mu d_{ij}) \right] X_{ij} - \sum_{j \in J} \lambda_j \tag{23}$$

subject to (11), (13), and (17).

We propose the following LB solution algorithm to efficiently solve the above LR problem and obtain the LB solution $\boldsymbol{X}_{LB}^k$ for each iteration *k*. As $\boldsymbol{X}_{LB}^k$ might be infeasible due to violating constraint (12), an upper-bound (UB) solution heuristic is proposed to fix the infeasibility and obtain a good feasible solution. Then the Lagrangian multiplier is updated by the subgradient method [32]. The proposed LR algorithm is labeled as Algorithm 1.

***LB solution algorithm***: Given a known $\boldsymbol{\lambda}^k$, we apply a greedy policy to assign vehicles to chargers one by one according to an increasing order with respect to the objective function value until all vehicles are assigned. To do so, a cost function $C(i,j)$ is defined as the cost of assigning vehicle $i$ to charger $j$ as in Eq. (24).

$$C(i,j) = t_{ij} + \theta_2 \max\left(t_j^A - t_{ij}, 0\right) + \frac{\theta_1}{\varphi_j}\left(e_i^* - e_i + \mu d_{ij}\right) \tag{24}$$

So the greedy policy assigns vehicle $i$ to the charger $j_i^*$ that minimizes the value of the objective function as in Eq. (25).

$$X_{ij_i^*} = 1, j_i^* = argmin_{j\in J_i}\left[C(i,j) + \lambda_j\right], \forall i \in I \tag{25}$$

where $J_i$ is the set of chargers that are reachable by vehicle $i$ given its current battery level, i.e., the subset of chargers satisfying constraint (13).

$$J_i = \left\{j \middle| e_i - \mu d_{ij} \geq e_{min}, \forall j \in J\right\}, \forall i \in I \tag{26}$$

The obtained LB solution $\boldsymbol{X}_{LB}(\boldsymbol{\lambda}^k)$ at iteration $k$ is the optimal solution of the LR problem. Note that for a P2J problem ($|J| < |I|$), a similar greedy policy applies by assigning chargers to vehicles until all chargers are assigned.

***UB solution heuristic***: We develop two distinguished heuristics to repair the feasibility of LB solutions and build UB solutions for the problems of P2 and P2J accordingly. The developed heuristics are described in Algorithm 2. Given an LB solution and constraint (12) for the P2 problem, the UB algorithm removes vehicles with higher energy levels (more flexible) from over-assigned chargers (chargers with more than one assigned vehicle) to a pool of unassigned vehicles. Then these unassigned vehicles are inserted to non-occupied chargers one by one based on their remaining energy levels (the vehicle with the least remaining energy (less flexible) is inserted first) using a greedy insertion policy. Afterwards, a local search procedure is applied to improve the incumbent feasible solution. A similar algorithm design logic is applied for the heuristic to find a UB solution for the P2J problem.

We test the proposed model on an illustrative example to show the model property in Appendix B. For large-scale problems, we generate 9 subsets of problems with up to 1000 vehicles and chargers. Our computational study in Appendix B shows that the proposed LR algorithm can obtain near-optimal solutions and suitable for large-scale real-time application for EV charging station assignment.

**Algorithm 1: Lagrangian relaxation algorithm**

1: Input: $\lambda^0 = 0$, $k = 0$, $Z_{UB} = INF$, $Z_{LB} = -INF$, $0 < \delta < 2$, $\bar{\varepsilon}$, and $k^{max}$.
2: Solve the LR problem by the LB solution algorithm and obtain the LB solution $X_{LB}^k$.
3: Update the LB: If $Z_{LB}^k > Z_{LB}$, set $Z_{LB} = Z_{LB}^k$
4: Repair the infeasible LB solution with the UB solution heuristic and obtain a feasible UB solution $X_{UB}^k$.
5: Update UB: If $Z_{UB}^k < Z_{UB}$, set $Z_{UB} = Z_{UB}^k$.
6: Update Lagrangian multipliers based on the subgradient method:
Compute the step size $t^k = \frac{\delta(Z_{UB}-Z_{LB})}{\sum_{j\in J}\left(\sum_{i\in I} X_{LB}^k - 1\right)^2}$ and update the multipliers as $\lambda^{k+1} = max\{\lambda^k + t^k(\sum_{i\in I} X_{LB}^k - 1), 0\}$.
7: Evaluate the optimality gap $\varepsilon^k = \frac{Z_{UB}-Z_{LB}}{Z_{UB}}$. If $\varepsilon^k \leq \bar{\varepsilon}$ or $k = k^{max}$, stop; otherwise k:=k+1 go to step 2.
8: Output: $X^* = X_{UB}$.

**Algorithm 2: Heuristics to find UB solutions**

**// Heuristic to find UB solutions for P2 ($|\boldsymbol{I}| \leq |\boldsymbol{J}|$).**
1: Given the current lower-bound solution $X_{LB}(\lambda)$, a set of vehicles $I$ and set of chargers $J$, initialize unassigned vehicle list $\bar{I} = \emptyset$ and temporary solution $X_{temp} = X_{LB}(\lambda)$.
2: //remove vehicles from over-assigned chargers
3: Find the list of chargers $J_1$ with more than one assigned vehicle.

4: **for** all chargers $j \in J_1$

5: Sort $e_i$ for all vehicles assigned to charger $j$ in descending order.

6: Remove *k–1* vehicles with the highest $e_i$ from charger $j$ to $\bar{I}$, where $k$ is the number of vehicles assigned to charger $j$. Update $X_{temp}$ accordingly.

7: **End for**

8: //Assign unassigned vehicles to unoccupied chargers

9: Sort $e_i$ for all vehicles $i \in \bar{I}$ in ascending order and obtain $\bar{I}_sorted$ list.

10: **for** all $i \in \bar{I}_sorted$

11: Assign $i$ to an unoccupied charger $j$ that has the minimum $C(i,j)$, then update $X_{temp}$ accordingly.

12: **End for**

13: // Local search for P2

14: For any two assigned vehicles $(i_1, i_2)$ with their current assigned chargers $(j_1, j_2)$, if $C(i_1, j_1) + C(i_2, j_2) - C(i_1, j_2) - C(i_2, j_1) > 0$ and both $(i_1, j_2)$ and $(i_2, j_1)$ satisfy Eq. (13), exchage their current assigned chargers. Update $X_{temp}$ accordingly.

15: **Output**: $X_{UB} = X_{temp}$

**// Heuristic to find UB solutions for P2J (|I|>|J|)**

16: Given the current lower-bound solution $X_{LB}(\lambda)$, set of vehicles $I$, and set of chargers $J$, initialize unassigned charger list $\bar{J} = \emptyset$ and temporary solution $X_{temp} = X_{LB}(\lambda)$.

17: //remove chargers from over-assigned vehicles

18: Find the list of vehicles $I_1$ with more than one assigned charger.

19: **for** all vehicles $i \in I_1$

20: Sort $C(i,j)$ for all assigned chargers in descending order.

21: Remove *k–1* chargers with the highest $C(i,j)$ to $\bar{J}$, where $k$ is the number of chargers assigned to vehicle $i$. Update $X_{temp}$ accordingly.

22: **End for**

23: // insert unassigned vehicles to available chargers

24: **if** there are unoccupied chargers, **then**

25: Sort $e_i$ for all non-assigned vehicles in ascending order and obtain $\bar{I}_sorted$ list.

26: **for** all $i \in \bar{I}_sorted$

27: Insert $i$ to an unoccupied charger $j$ with minimum $C(i,j)$ if the assignment $(i,j)$ satisfies Eq. (13). Update unoccupied charger list and $X_{temp}$ accordingly. If each charger is assigned by one vehicle, break.

28: **End for**

29: **End if**

30: // repair infeasible solution if the above procedure fails

31: **while** the number of unoccupied chargers >0

32: Sort $e_i$ for all assigned vehicles in descending order.

33: Pop up the vehicle in the sorted assigned vehicle list with highest $e_i$, update unassigned vehicle and charger lists.

34: **for** all unassigned vehicles $i$ and unoccupied chargers $j$

35: Assign $i$ to $j$ with the minimum $C(i,j)$ while satisfying Eq. (13). Update the unassigned vehicle and charger lists and $X_{temp}$ accordingly. If each charger is occupied by one vehicle, leave the loop.

**End for**

36: **End while**

37: // local search

38: For any unassigned vehicles, exchange with current assigned vehicles if the resulting solution decreases the objective function value.

39: Apply the exchange procedure of step 14 to improve the current solution.

40: **Output**: $X_{UB} = X_{temp}$.

# Computational study for dynamic dial-a-ride services using EVs in Luxembourg

The proposed methodology is applied for a realistic dynamic dial-a-ride service using EVs in Luxembourg. The goal is to demonstrate the benefit of the proposed methodology for reducing the total charging operation delay and assess its impact on system performance. The computational study is implemented on the simulation platform previously used for dynamic ridesharing with transit transfers [33, 34] but is extended to handle EVs recharging as a multi-server queuing system.

# Dynamic dial-a-ride simulation platform using EVs

Consider a TNC operating a fleet of homogeneous EVs to provide dial-a-ride service in Luxembourg. Ride requests are unknown in advance and arrive stochastically. The operator makes vehicle dispatch and routing decisions over time. We adopt a non-myopic vehicle dispatching policy to anticipate future system delays [34]. This policy minimizes the marginal system cost increases when inserting a new customer on an existing tour by considering the future system cost as an M/M/1 queue delay. A traveling salesman problem with pickup and delivery problem is solved by a re-optimization-based heuristic to insert a new request on existing routes which satisfies vehicle capacity constraint and precedence constraint (a pick-up location is visited before its corresponding drop-off location). Following the previous studies [34], no time window constraints are associated with customer requests. The applied vehicle routing policy can be substituted by other approaches considering time windows and other vehicle operational constraints. In an EV-enabled dial-a-ride system, the vehicle's remaining energy level needs to be considered when assigning a new request to an existing vehicle tour. We assume that the remaining energy of a vehicle after serving all customers (including the new request) and returning to its depot needs to be no less than the minimum reserve energy $e_{min}$, as in Eq. (27):

$$e(v,\bar{x}_t^v) - e(v,x_t^v) \geq e_{min} \tag{27}$$

where $x_t^v$ is the current tour of vehicle $v$. $\bar{x}_t^v$ is the post-evaluated tour after inserting the new request. When dispatching vehicles to pick up a new customer, the dispatching center first determines a list of energy-feasible vehicle candidates (satisfying Eq. (27)). If no vehicles are energy-feasible, the new request is rejected. Otherwise, the customer is assigned to the vehicle with lowest marginal system cost, as in Eq. (28).

$$\{v^*, x_t^{v*}\} = \text{argmin}_{v\in V',x}[c\ (v,\bar{x}_t^v) - c(v,x_t^v\ )] \tag{28}$$

where $V'$ is the set of vehicles satisfying Eq. (27). $c(v,x)$ is a cost function with service tour $x$ defined as Eq. (29).

$$c(v,x) = \gamma T(v,x) + (1-\gamma)\left[\beta\ T(v,x)^2 + \sum_{n\in P_v} \bar{Y}_n(v\ ,x)\right] \tag{29}$$

where $T(v,x)$ is the travel time of tour $x$. $\bar{Y}_n(v\ ,x)$ is the journey time (i.e., waiting time plus in-vehicle travel time) for passenger $n$ among the set of passengers $P_v$ assigned to vehicle $v$. Parameter $\alpha$ is a weight considering the trade-off between operational cost and customer inconvenience. $\beta$ is a parameter between 0 and 1. When $\beta = 0$, the resulting vehicle dispatching policy is myopic since it does not consider future approximate system delays for the current vehicle dispatching decision.

For the recharging policy, the entire planning period is discretized into a set of charging decision epochs with the time interval $\Delta$. The latter is set as $\frac{e_{max}-e_{min}}{\varphi_{max}}$ to allow EVs to be recharged to the desired energy level ($e_{max}$) within a charging decision epoch. Note that the impact of epoch length will be analyzed in Section 4.3.4. $\varphi_{max}$ is the maximum charging rate of all chargers in the studied area.

We first solve the P1 problem for each vehicle and obtain its optimal charging plan for each decision epoch. The inputs of the P1 problem are the probability distributions of vehicles in the driving state and expected waiting times at charging stations, which can be obtained from the historical driving and charging patterns of vehicles. These data can be easily obtained by the operator using dedicated fleet management software with GPS tracking. The output of the P1 problem provides the vehicle-specific charging schedule (plan). Then we solve the P2/P2J problem for each charging decision epoch. Due to the stochastic nature of customer arrivals and charging station occupancy, a vehicle's exact charging time and amount of recharged energy depends on the charging rate of its assigned charging station and the battery level when arriving at a charging station. Moreover, vehicles can go to recharge only after all customers on board are served. We apply the discrete event simulation technique to include queueing delays at charging stations.

# Luxembourg case study

We apply the proposed methodology to a dynamic dial-a-ride case study using EVs in Luxembourg. Such services using conventional gasoline shuttles have been operating in Luxembourg (e.g. Flexibus, https://www.sales-lentz.lu/en/individuals/shuttle-upon-request/). Luxembourg has promoted a sustainable mobility initiative aimed at shifting mobility practices from the current high car dependency towards multimodal and soft mobility alternatives. To promote e-mobility, the government plans to install 1600 charging points, named Chargy (https://chargy.lu/), with 22 kWh Level 2 chargers in the entire country in order to meet the future charging needs of electric/hybrid vehicles. Currently, a total of 814 Level 2 charging plugs have been installed. We consider a TNC that provides dynamic dial-a-ride services in Luxembourg using EVs and the Chargy network for recharge in the daytime. For simplicity, the use of charging infrastructure from other EVs (private EVs or other types of EVs) is not considered in this study. The impact of the charging needs of the other stochastically arriving EVs will be discussed in Section 4.4.

We use the Luxmobil open trip dataset obtained from the 2017 Luxmobil survey [35] (for which 40,000 households in Luxembourg and 45,000 cross-border workers were surveyed in 2017, with response rates ranging from 26% to 30%) to generate customer demand in Luxembourg. The trip dataset contains 82000 trips representing typical one-workday trips for Luxembourg residents and its cross border workers from France, Germany, and Belgium. The trip dataset is anonymized and provides trip information related to departure time, origin, and destination at the level of neighborhood, transportation means, etc. The extrapolation from the Luxmobil survey data to the entire population (i.e. Luxembourg residents and its cross border workers) was done for German, French and Belgian commuters by the number of workers in Luxembourg per municipality/district. Residents of Luxembourg are extrapolated based on the number of inhabitants and the age distribution of each municipality/district. We randomly generate 1000 customer ride requests from 6:30–22:00 in Luxembourg. The impact of different demand intensities is evaluated in the Sect. Sensitivity analysis. As the survey data contains trip origins and destinations at the neighborhood level, we use another geo-referenced data set (https://data.public.lu/fr/datasets/adresses-georeferencees-bd-adresses/) to randomly generate the geographical coordinates within the same neighborhoods. The geo-referenced data set includes all neighborhoods in Luxembourg and contains all building addresses appearing in the national register of streets and neighborhoods. Fig 2 shows the heat map of a one-day data set of 1000 customers' pickup locations in the study area. We can observe that most trips are distributed between Luxembourg City and the cities of Esch-sur-Alzette, Dudelange, and Mersch. Fig 3 depicts the distribution of customer arrival times. We observe that a morning arrival peak occurs between 7:00 and 8:00, and an afternoon peak occurs between 17:00 and 19:00.

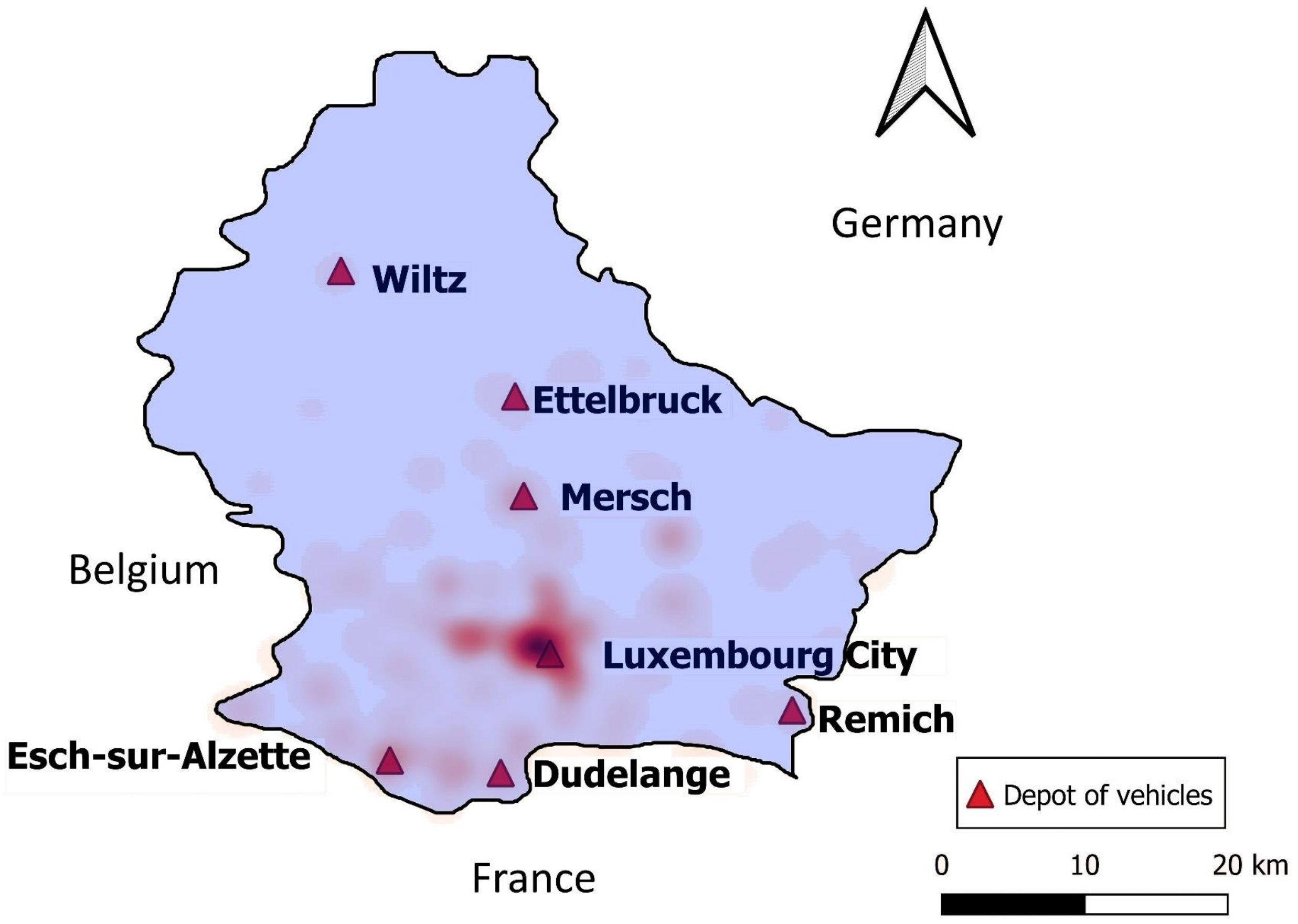

**Fig 2. Heat maps of customer pickup location between 6:30–22:00 in the study area.**

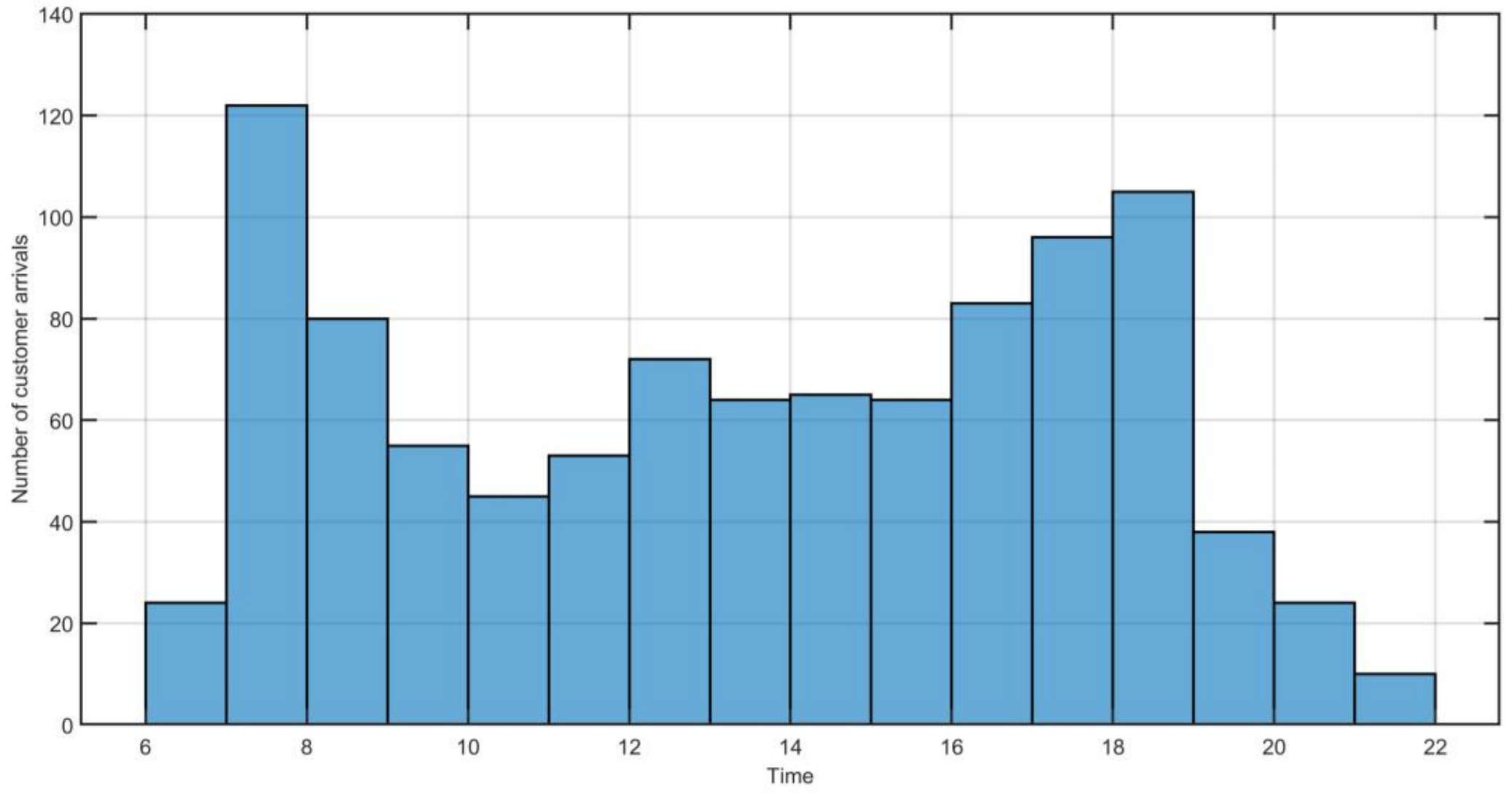

**Fig 3. Example of 1000 customer arrivals from 6:30–22:00 in the study area.**

The fleet size is assumed to include 50 homogeneous 8-seater fully electric shuttles, representing 20 customers/vehicle/day. The characteristics of the EVs are based on Volkswagen's 8-seat 100% electric Tribus (https://www.tribus-group.com/zero-emission-volkswagen-e-crafter-electric-wheelchair-minibus/). The battery range of the Tribus is 35.8 kWh, with a practical range up to 150 km. The average energy consumption per kilometer travelled is assumed to be constant. For practical applications, a drive cycle can be applied instead of a constant energy consumption policy. We assume that EVs are recharged to 80% of their battery capacity to maximize battery lifetimes. Vehicles are initially located at 7 different depots around the municipality centers of Luxembourg City, Esch-sur-Alzette, Ettelbruck, Dudelange, Mersch, Remich, and Wiltz. Vehicle dispatching and routing policy is based on the non-myopic policy in Section 4.1.

On the charging infrastructure side, we test the proposed methodology based on three different charging infrastructure configuration scenarios as follows.

- Scenario 1 (Chargy only): EVs can only be recharged in the daytime on Luxembourg's current public Chargy infrastructure, which comprises 814 level 2 (L2) chargers with 22 kWh power. Based on this scenario, recharging EVs using an L2 charger from 0% to 80% require around 1.3 hours. Fig 4 shows the spatial distribution of the charging stations in Luxembourg.
- Scenario 2 (DC fast only): EVs recharge exclusively on 9 DC fast chargers with 50 kWh power, distributed in different municipalities in Luxembourg (see Fig 4).
- Scenario 3 (Chargy+DC fast): EVs can be recharged both on the Chargy network or at DC fast charging stations.

Note that we assume that all charging stations are available at the beginning of each day. EVs are initially charged to 80% at their depot and maintain an energy level no less than 10% for customer service operations. For practical applications, higher minimum charge levels can be applied if necessary.

The charging decision epoch is 30 minutes, reflecting the charging time for a Tribus vehicle from $e_{min}$ (10% of battery capacity) to $e_{max}$ (80% of battery capacity) using a 50 kWh DC fast charger. We randomly generate 13 independent 1000-customer demand data sets from the Luxmobil open trip dataset. A total of 10 simulation runs on distinguished demand data are conducted to derive the probability distribution of a vehicle being in the driving state and the expected waiting time when arriving at charging stations using a need-based charging policy (reference policy). In practice, such information can be easily collected by an operator based on vehicles' historical driving patterns and the waiting times experienced at charging stations. The **reference policy** states that EVs go to the nearest unoccupied charger to recharge their batteries to 80% whenever the battery level is lower than 20% after serving all customers on board [16, 21]. The three remaining demand data sets are used to test the performance of the proposed charging strategy.

Table 2 shows the details of the simulation parameter settings and the parameters used for the charging schedule and charging station assignment models.

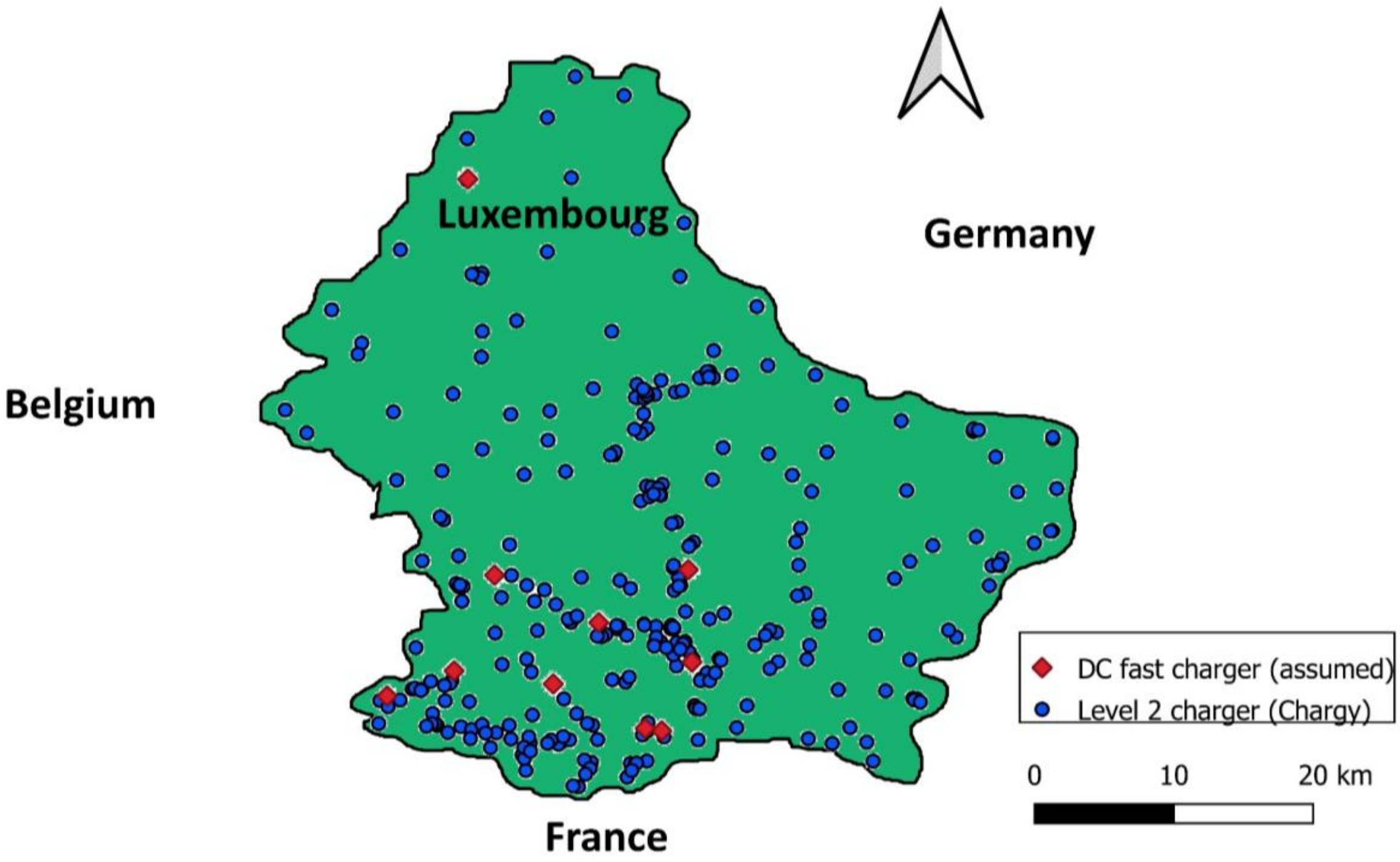

**Fig 4. Charging station locations in the study area.**

**Table 2. List of parameters used for the Luxembourg case study.**

| | | | |
|---|---|---|---|
| Number of customers | 1000 | $\rho$ | 0.2485 (euro/min.)[2] |
| Number of vehicle depots | 7 | $\vartheta$ | 0.2756 (euros/kWh)[3] |
| Fleet size | 50 | $\mu$ | 0.2387 (kWh/km)[4] |
| Capacity of vehicles | 8 pers./veh. | $\bar{c}$ | 5.77 euros[5] |
| Vehicle speed | 50 km/hour | $\Delta$ | 30 min. |
| Battery capacity | 35.8 kWh[1] | $T$ | 6:30–22:00 |
| Number of chargers | | $\varphi_{L2}$ | 22/60 (kW/min.) |
| Level 2 (22 kWh) | 814 | $\varphi_{\text{DC fast}}$ | 50/60 (kW/min.) |
| DC fast (50 kWh) | 9 | | |
| Battery capacity | 35.8 kWh | $\beta$ | 0.025 |
| $e_{min}$ | $0.1B$ | $\gamma$ | 0.5 |
| $e_{max}$ | $0.8B$ | | |

Remark: 1. EV characteristics are based on a Volkswagen-powered 100% electric minibus with a 150 km range (https://www.tribus-group.com/zero-emission-volkswagen-e-crafter-electric-wheelchair-minibus/). 2. Equivalent to a 14.9 euros/hour wage rate. As a reference, the range of the gross salary of bus drivers in Luxembourg is 2474–4846 euros/month (around 14.4–28.2 euros/hour, https://www.paylab.com/lu/salaryinfo). 3. Electricity price is based on the current fare using Chargy's plugs (https://www.eida.lu/en/chargy) with VAT. 4. Based on the driving efficiency of the Tribus ($\mu = 35.8/150$). 5. Based on the energy consumption cost of the average access distance to charging stations of 5 km.

# Results

We first derive the charging plan for each vehicle under different charging infrastructure scenarios. Then we run the dynamic dial-a-ride simulation under the different charging policies and charging infrastructure configuration scenarios. The discrete probability distribution is estimated using the frequency of a vehicle being in the driving state over 10 runs using independent random demand samples as previously mentioned. The charging policy is based on the reference policy using the Chargy network only for vehicle recharge. Fig 5 shows an example of the empirical probability distribution that a vehicle is in the driving state. Here we select two vehicles with quite different driving patterns. Vehicle 1 is initially located in Luxembourg City, with higher customer demand, while vehicle 40 with its depot at Mersch has lower customer demand.

Fig 6 reports the expected waiting time to be served when arriving at charging stations in each charging decision epoch. We vary the charging infrastructure configuration scenarios to obtain the respective expected waiting times under each scenario. We can observe that when using the Chargy network only, a peak in the expected waiting time of 30 minutes appears around at 9:00, which reduces to less than 10 minutes after around 10:30. When adding 9 DC fast chargers to the existing Chargy network, the peak in the expected waiting time at 9:00 is reduced to around 15 minutes. We find that when all vehicles recharge only at DC fast chargers, the profile of the expected waiting time in each charging decision epoch is higher than the two other scenarios. Fig 7 shows an example of the optimal charging plan obtained from P1 for vehicle 1 and vehicle 40. We can observe that vehicle 1 needs to be recharged 4 times to nearly 80%, at around 10:00, 13:30, 16:30, and 19:30, whereas vehicle 40 requires 3 recharges, at around 10:30 and 15:30 to around 80% and at 22:30 to around 15%.

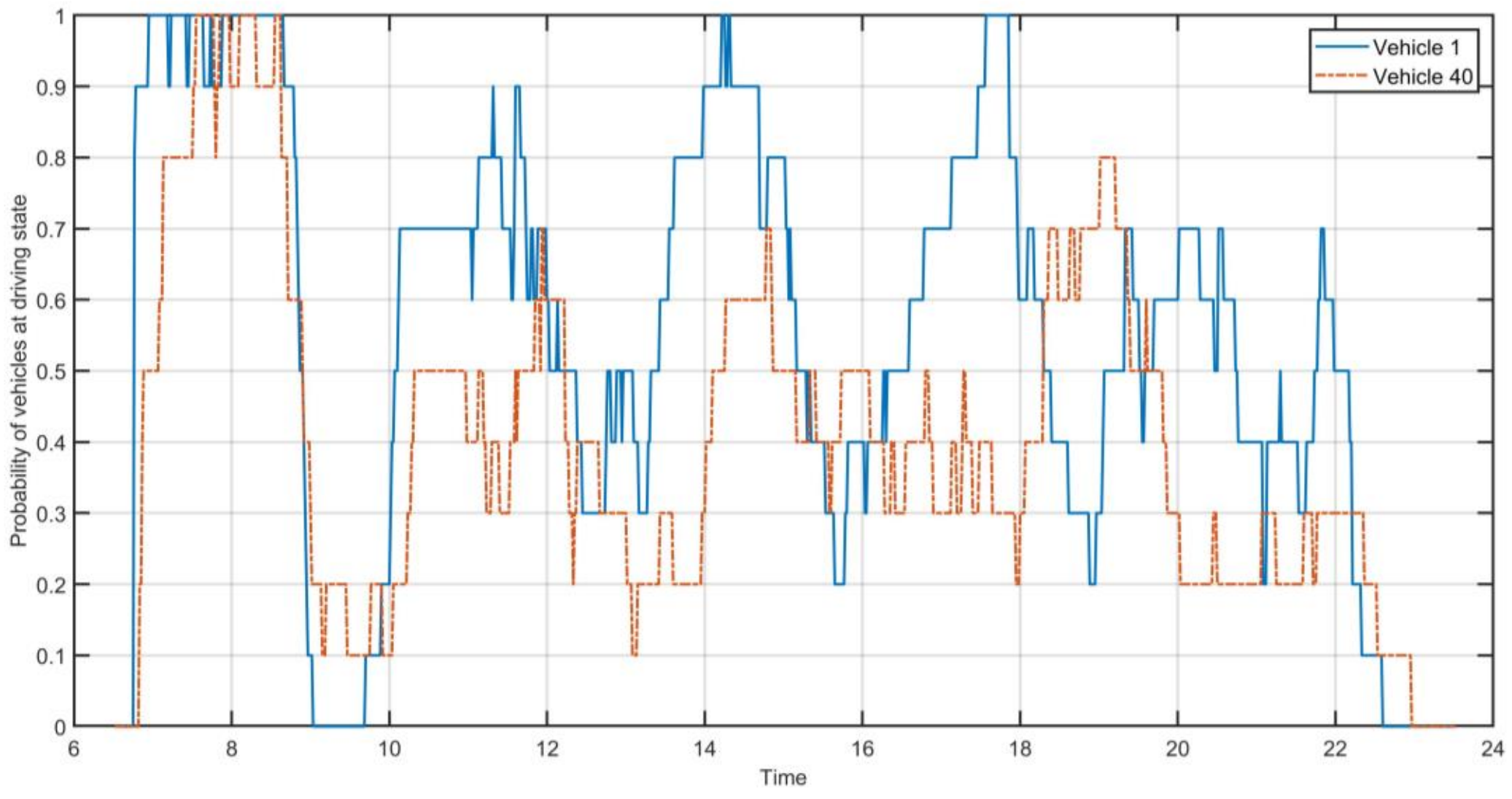

**Fig 5. Example of the probability distribution of a vehicle being in the driving state.**

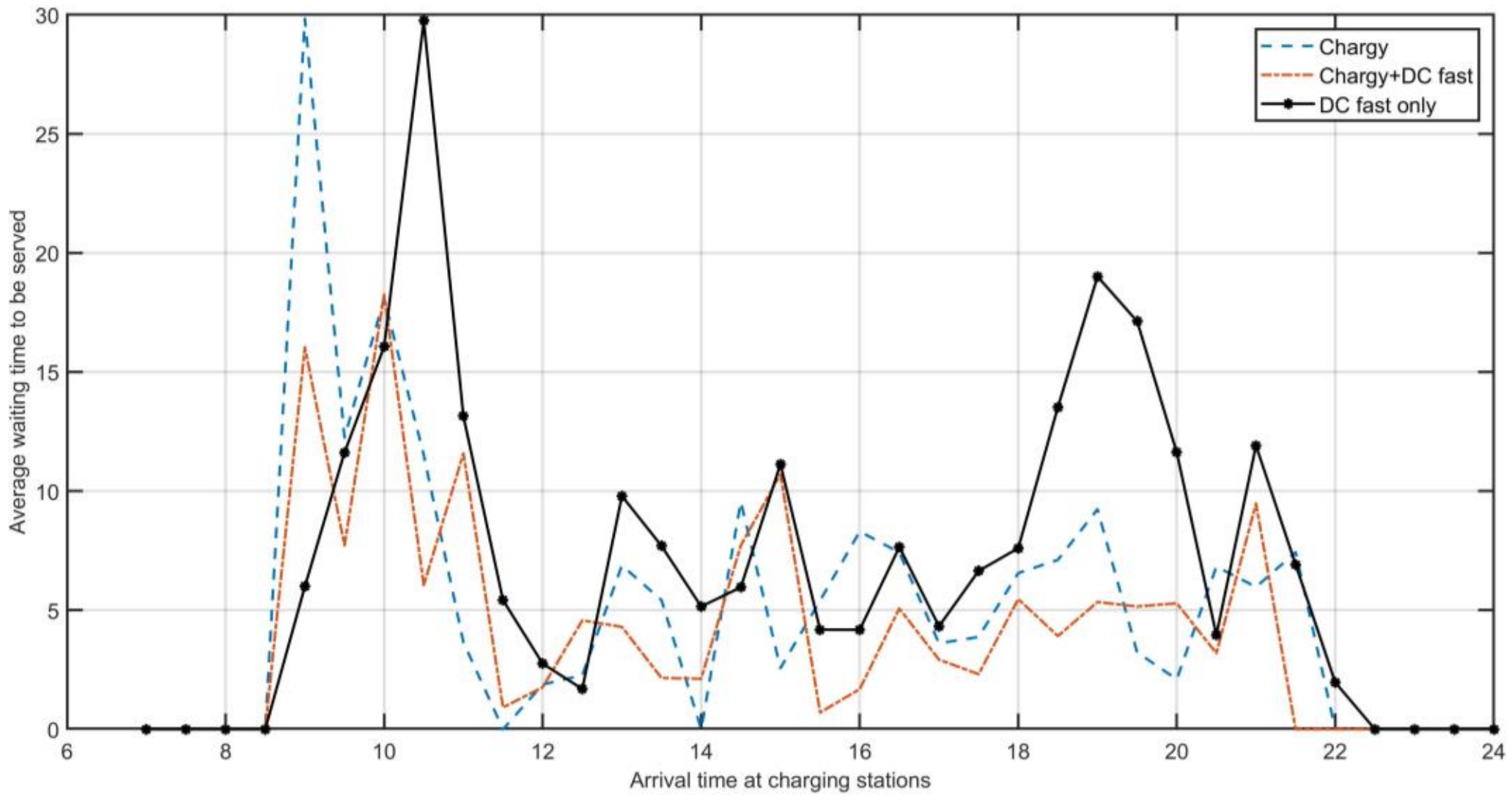

**Fig 6. Expected waiting times to be served when arriving at charging stations under different charging infrastructure scenarios.**

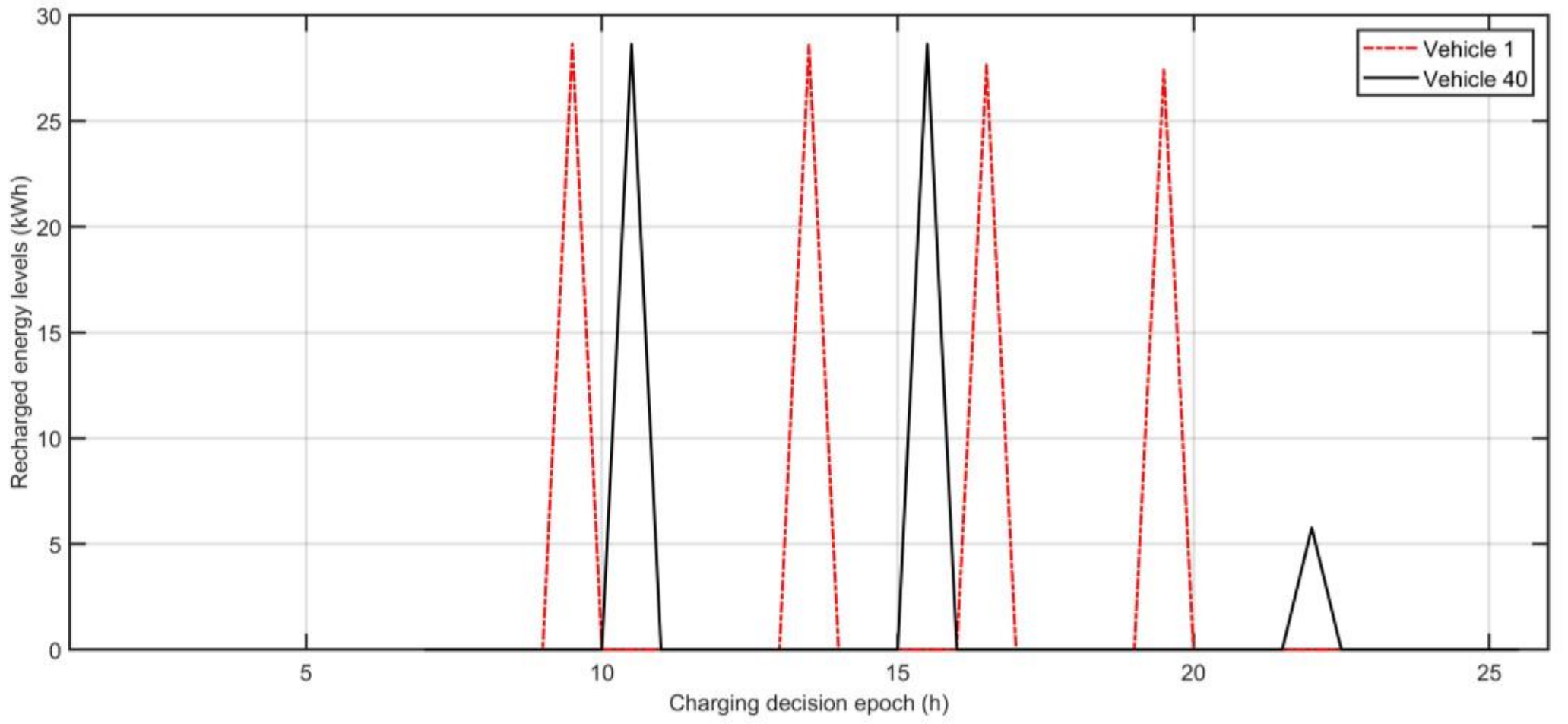

**Fig 7. Example of vehicle optimal charging plans in different charging decision epochs.**

To assess the proposed charging policy, two charging reference policies are considered as follows.

- **Need-based nearest charging station assignment policy (NS)**: vehicles recharge to 80% at the nearest unoccupied charger whenever their battery range is lower than 20% after serving onboard customers.
- **First-come-first-served (FCFS) minimum charging delay policy**: vehicles go to recharge whenever their battery range is lower than 20% after serving onboard customers. The assigned charger $j$ for vehicle $i$ is based on the lowest estimated charging operation time at time $t$ as in (21).

$$j^*(i,t) = argmin_{j\in J}\{t_{ij} + \widetilde{W}_{ij}(t)\}, \tag{30}$$

where $t_{ij}$ is the travel time from the location of vehicle $i$ to charger $j$. $\widetilde{W}_{ij}(t)$ is the expected waiting time estimated as the difference between the earliest available time of chargers and the arrival time of vehicle $i$.

We refer to the proposed optimal charging scheduling and assignment policy as the **OCP policy**. The computational results are based on the average performance using three independent demand data sets**.**

Table 3 shows the impact of the three charging policies on the total charging delays and the amount of recharged energy under different charging infrastructure scenarios. The first three columns report the average charging waiting time, charging time, and vehicle idle time for recharge per vehicle charging. The other columns report the measures related to the overall charging costs/times for the fleet for a full day of operation. First, for each vehicle charging operation, the OCP policy shows the best performance, with the least waiting time, charging time, and vehicle idle time for recharge over all three scenarios. The average charging waiting time over the three scenarios is 3.1 min., compared to the 9.9 min. (–68.9%) of NS and 11.3 min. for FCFS (–72.6%). Regarding the average charging time, the OCP policy results in 37.4 min., compared to NS with 51 min. (–26.8%) and FCFS with 52.4 min. (–28.6%). Similar results can be found for the average vehicle idle time per vehicle recharge, with an average savings of 29%–32.7%.

Second, in terms of the total waiting time of the fleet, the results show that using the NS policy would lead to a high charging delay, with 27.2 hours for Chargy only, 24.3 hours for DC fast only, and 21.7 hours for Chargy plus DC fast. Adopting the FCFS policy would not significantly reduce the waiting time compared to the NS policy. However, the benefit of the OCP policy is very significant in terms of reducing the charging operation delay: on average, –73.5% compared to that of the NS policy and –76.5% compared to the FCFS policy over the three scenarios. In terms of the total charging time of the fleet, it can be observed that using Chargy only, the NS policy would lead to a high value of 159.6 charging hours, compared to that of using DC fast chargers only (69 hours) and Chargy+DC fast chargers (148.7 hours). Using the OCP policy would lead to significant total charging time savings under the three charging infrastructure scenarios, and in particular the Chargy+DC fast charger scenario. On average, the benefit in terms of reducing the total charging time is –38.1% compared to the NS policy and –39.0% compared to the FCFS policy. Finally, using the OCP policy can significantly reduce the amount of total charged energy and energy costs. On average, the total charged energy cost saving is around 27% compared to the other two policies.

**Table 3. Impact of different charging policies on total system delays, recharged energy, and costs.**

| Scenario | Charging policy | Average charging waiting time[1] (min.) | Average charging time[1] (min.) | Average operational time for charge[2] (min.) | Total waiting time[3] (hour) | Total charging time[3] (hour) | Total amount of charged energy[3] (kWh) | Total charged energy cost[3] (euro) |
|---|---|---|---|---|---|---|---|---|
| Chargy only | NS | 11.1 | 65.0 | 76.4 | 27.2 | 159.6 | 3512.1 | 967.9 |
| | FCFS | 9.7 | 65.0 | 75.5 | 23.5 | 158.1 | 3477.3 | 958.3 |
| | OCP | 4.7 | 53.6 | 61.9 | 9.9 | 112.5 | 2475.9 | 682.3 |
| DC fast only | NS | 10.0 | 28.4 | 45.3 | 24.3 | 69.0 | 3452.2 | 951.4 |
| | FCFS | 17.0 | 28.5 | 54.6 | 41.8 | 70.1 | 3504.1 | 965.7 |
| | OCP | 3.7 | 24.8 | 35.4 | 7.8 | 52.3 | 2614.2 | 720.5 |
| Chargy+DC fast | NS | 8.7 | 59.7 | 69.4 | 21.7 | 148.7 | 3567.4 | 983.2 |
| | FCFS | 7.1 | 63.6 | 71.6 | 17.3 | 155.1 | 3493.2 | 962.7 |
| | OCP | 0.8 | 33.7 | 38.4 | 1.7 | 69.0 | 2534.0 | 698.4 |
| Average over three scenarios | NS | 9.9 | 51.0 | 63.7 | 24.4 | 125.8 | 3510.6 | 967.5 |
| | FCFS | 11.3 | 52.4 | 67.2 | 27.5 | 127.8 | 3491.5 | 962.3 |
| | OCP | 3.1 | 37.4 | 45.2 | 6.5 | 77.9 | 2541.4 | 700.4 |

Remark: [1] Measured as the time per recharge per vehicle. [2] Includes vehicle travel time to reach a charger, charging waiting time, and charging time. [3] Measured for the e-fleet for a full day of operation.

Table 4 shows the system performance from the customer perspective under different charging policies. It can be observed that using the NS and FCFS policies, almost all customers (99.5%) are served given different charging infrastructure scenarios. The OCP policy has, on average, a 94.1% (–5.5%) customer service rate. This is due to the stochasticity of ride requests; the charging plans of the OCP policy obtained from historical vehicle driving patterns might not be able to fit certain long-trip requests perfectly, resulting in rejections due to insufficient energy remaining in the vehicles. This trade-off between the service quality and efficiency can be improved by increasing the fleet size or using a mixed fleet of gasoline and electric vehicles. Second, the results indicate that for the DC-fast-only scenario, customer inconvenience can be improved due to the overall charging time savings given the current charging demand of the fleet. In terms of customer inconvenience, adopting the OCP policy would increase customer waiting times by 2.7 minutes on average over three scenarios and passenger journey times by 5 minutes compared to the NS and FCFS policies.

We can conclude that adopting the OCP policy can lead to a significant reduction in charging operation delays and costs compared to the NS and FCFS policies, while maintaining a high customer service rate with very limited perturbation in terms of customer inconvenience.

**Table 4. System performance under different charging policies.**

| Scenario | Charging policy | Mean passenger waiting time | Mean passenger journey time | Mean vehicle travel time | % of customers served |
|---|---|---|---|---|---|
| Chargy only | NS | 13.3 | 34.0 | 393.5 | 99.3 |
| | FCFS | 13.1 | 33.8 | 392.5 | 99.3 |
| | OCP | 15.8 | 39.0 | 367.7 | 92.9 |
| DC fast only | NS | 10.0 | 29.3 | 392.2 | 100.0 |
| | FCFS | 10.6 | 30.4 | 390.9 | 100.0 |
| | OCP | 14.1 | 36.6 | 380.9 | 96.2 |
| Chargy+DC fast | NS | 13.4 | 34.2 | 396.8 | 99.5 |
| | FCFS | 13.0 | 33.9 | 392.1 | 99.3 |
| | OCP | 14.9 | 37.4 | 374.7 | 93.2 |
| Average over three scenarios | NS | 12.2 | 32.5 | 394.2 | 99.6 |
| | FCFS | 12.2 | 32.7 | 391.8 | 99.5 |
| | OCP | 14.9 | 37.6 | 374.4 | 94.1 |

NS: nearest station charging policy; OCP: optimal charging plan policy; AOCP: adapted optimal charging plan policy. Passenger journey time includes passenger waiting time and in-vehicle travel time.

# Sensitivity analysis

To investigate the impact of model parameters on the performance of the proposed method, a sensitivity analysis is designed to answer the following questions.

- What is the impact of the length of charging decision epochs on the charging queueing delay and customer inconvenience?
- How do ride demand changes affect the total charging delays and costs?
- What are the benefits of increasing the vehicle battery range in terms of reductions in charging delays, costs, and customer inconvenience?

The reported results are based on the average of 3 runs using the three independent customer demand test data sets on the Chargy+DC fast chargers scenario.

## a) Impact of the length of the charging decision epochs

The length of the charging decision epoch determines the frequency of vehicles scheduled for recharge, which may influence the effectiveness of vehicle charge scheduling and delays. A longer decision epoch would postpone the charging operations of vehicles, while a shorter decision epoch would limit the amount of energy to be charged within one epoch and increase the number of charging operations and the vehicle idle time due to charging operations. We vary the length of charging decision epochs from 10 to 60 minutes with a 10-minute interval and assess the impact on the total charging delay and energy costs. The results are shown in Table 5 and Fig 8. Two insights can be drawn: 1) The length of the charging decision epoch impacts the total charging delay and energy costs. When the epoch length is too short (less than 30 minutes), the maximum amount energy of energy that can be charged within a decision epoch is constrained, compromising the obtained charging plans of vehicles. Consequently, the vehicles need to charge more frequently, resulting in higher charging delays and a lower customer service rate; 2) when the decision epoch is greater than a critical value (30 minutes, allowing vehicles to be charged to 80% on a DC fast charger), the charging delays and energy costs are similar thanks to the adapted vehicle driving probability and charging waiting time in each epoch. This suggests that the charging decision epoch length should be long enough to allow vehicles to be charged to the desired energy level within one epoch using the faster chargers.

**Table 5. The impact of the charging decision epoch length on the performance of the OCP policy.**

| Charging decision epoch length | Total waiting time (hour) | Total charging time (hour) | Total amount of energy charged (kWh) | Total charged energy cost (euro) | Mean passenger waiting time (min.) | Mean passenger journey time (min.) | % of customers served |
|---|---|---|---|---|---|---|---|
| 10 | 5.9 | 23.7 | 932.1 | 256.9 | 17.1 | 38.5 | 54.0 |
| 20 | 5.3 | 66.0 | 2300.5 | 634.0 | 15.9 | 38.2 | 86.9 |
| 30 | 1.7 | 69.0 | 2534.0 | 698.4 | 14.9 | 37.4 | 93.2 |
| 40 | 1.6 | 73.1 | 2485.6 | 685.0 | 15.1 | 37.3 | 92.6 |
| 50 | 1.3 | 69.4 | 2548.7 | 702.4 | 14.6 | 36.6 | 94.7 |
| 60 | 1.2 | 74.5 | 2461.7 | 678.4 | 14.3 | 37.1 | 92.8 |

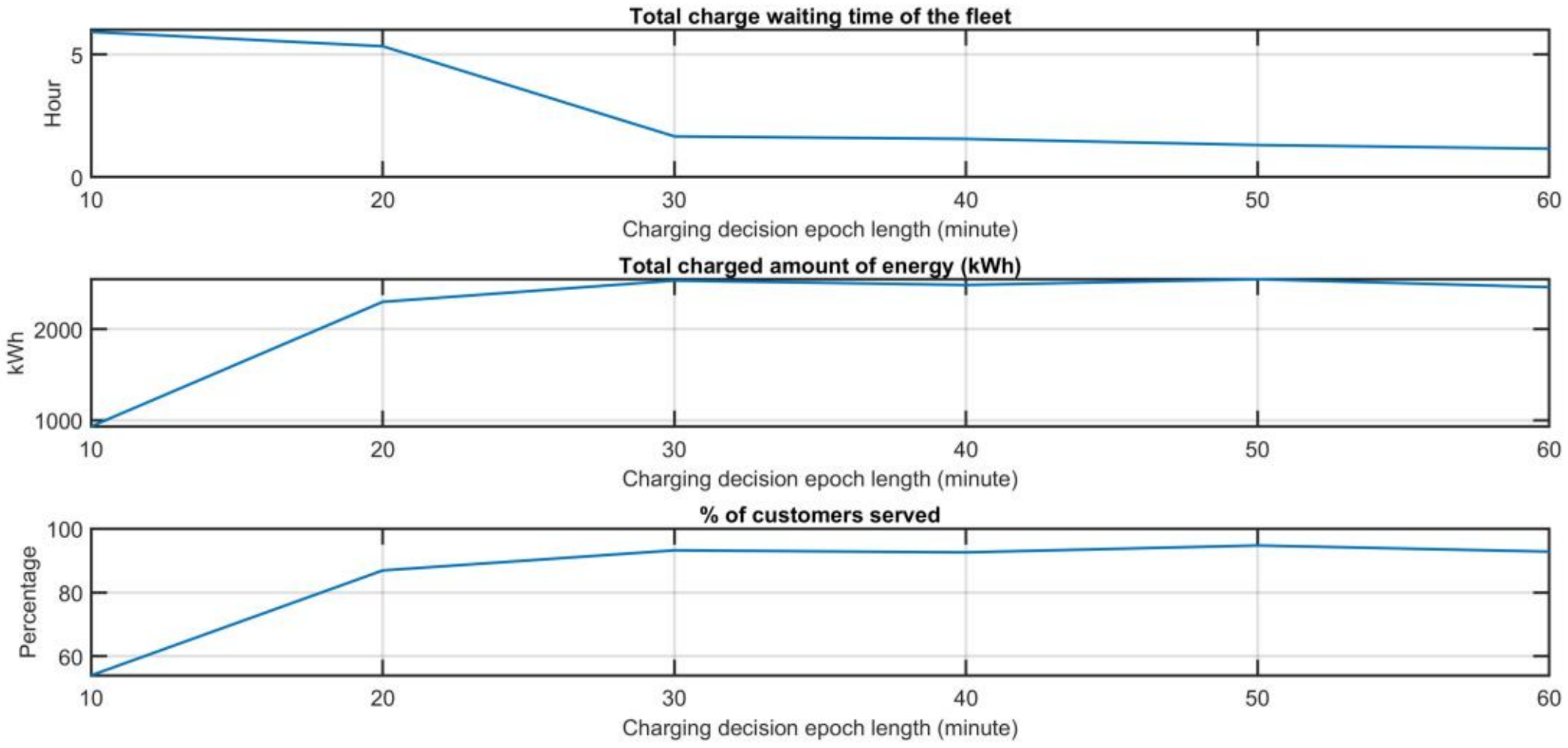

**Fig 8. The impact of charging decision epoch length on the total charging waiting time of the fleet, the total amount of charged energy, and the percentage of served customers.**

## b) Impact of customer demand variation

We vary customer demand as 500, 1000, 1500, and 2000 trip requests, respectively, randomly sampled from the Luxmobil open trip dataset. The results in Table 6 indicate that using the OCP policy can significantly reduce the total charging waiting times of the fleet compared to the NS and FCFS policies, given different customer demand. We found that when demand increases, the total amount of energy charged may vary depending on the total distance traveled by vehicles to serve customers. The OCP policy allows reducing the total energy costs of vehicle charges from –12.0% to –38.6% compared to the NS and FCFS policies, depending on the different customer demand intensities and service rates.

**Table 6. Total charging delays and costs of the fleet given different charging policies and customer demand.**

| Demand | 500 | 1000 | 1500 | 2000 | 500 | 1000 | 1500 | 2000 | 500 | 1000 | 1500 | 2000 |
|---|---|---|---|---|---|---|---|---|---|---|---|---|
| Charging policy | Total charging waiting time (hours) | | | | Total charged energy cost (euro) | | | | % of customers served | | | |
| NS | 4.1 | 21.7 | 8.3 | 8.4 | 493.4 | 983.2 | 1039.5 | 739.9 | 100.0 | 99.5 | 78.9 | 51.2 |
| FCFS | 5.9 | 17.3 | 12.1 | 14.4 | 483.9 | 962.7 | 1124.8 | 768.5 | 100.0 | 99.3 | 79.5 | 51.0 |
| OCP | 1.0 | 1.7 | 5.5 | 1.2 | 302.8 | 698.4 | 880.3 | 651.3 | 96.8 | 93.2 | 77.9 | 50.4 |

## c) Impact of battery range

To evaluate the impact of the vehicle battery range, we extend the battery capacity of the current Tribus from 35.8 kWh to 53.7 kWh (+50%) and 71.6 kWh (+100%), respectively. The results in Table 7 show that extending the battery range would significantly reduce the total waiting times of the NS policy and FCFS policy, while applying the OCP policy would increase the total waiting times from 1.7 hours to 3.6 hours due to the longer charging times of vehicles. The gains from the reduction of total charged energy (costs) are very significant for the OCP policy (from 2534 kWh to 1372 kWh, –45.8%). Moreover, the customer service rate is improved from 93% to 98%, close to the other charging policies. However, the NS and FCFS policies still charge a similar amount of energy (3352 kWh on average) over different battery capacities due to the 80%-charge policy for vehicle charging operations.

**Table 7. The impact of battery range on the performance of the OCP policy.**

| Battery capacity (KWh) / Charging policy | 35.8 | 53.7 | 71.6 | 35.8 | 53.7 | 71.6 | 35.8 | 53.7 | 71.6 |
|---|---|---|---|---|---|---|---|---|---|
| | Total charging waiting time (hours) | | | Total amount of energy charged (kWh) | | | Customer service rate (%) | | |
| NS | 21.7 | 9.3 | 5.1 | 3567.4 | 3364.0 | 3145.4 | 100 | 100 | 100 |
| FCFS | 17.3 | 17.3 | 7.9 | 3493.2 | 3320.5 | 3219.8 | 99 | 100 | 100 |
| OCP | 1.7 | 3.0 | 3.6 | 2534.0 | 1913.1 | 1372.3 | 93 | 97 | 98 |

# Discussion

The computational study first illustrates the characteristics of the charging scheduling model and the vehicle–charger assignment on a small example. Then the LR algorithm for solving online vehicle–charger assignment is evaluated on several numerical test instances. A realistic dynamic dial-a-ride service case study in Luxembourg is designed to assess the performance of the proposed approach and compare it with two widely used charging policies. A number of insights can be summarized as follows.

- The single-vehicle charging scheduling problem for dynamic shared on-demand mobility services can be decomposed into a multi-stage vehicle battery recharge problem to determine when and how much energy to charge in each charging decision epoch. The objective is to minimize charging operational costs while meeting vehicle driving needs for the next stage and battery-level-related constraints. The minimum amount of energy to charge for each vehicle is solved sequentially to obtain vehicle daily charging schedules that consider the expected queuing delays for charging and stochastic driving needs over time.
- The real-time vehicle–charger assignment model considering the current charging system queuing states considerably reduces vehicle charging waiting times and vehicle idle times for recharge in a dynamic environment. The LR algorithm allows for solving the mixed-integer assignment problem for large-scale test instances with 1000 vehicles and 1000 chargers within 3 minutes, with an optimality gap of 0.5%. The algorithm is suitable for the real-time vehicle-charger assignment of electric fleet charging operations to minimize total vehicle charging delays.
- The realistic dynamic dial-a-ride case study in Luxembourg under different charging infrastructure settings shows that the proposed charging schedule policy can reduce, on average, the total charging waiting time (–74.9%), charging time (–38.6%), and charging cost (–27.4%) compared to the nearest charging station charging policy and the minimum charging delay policy.
- A sensitivity analysis provides insight into the impact of the length of the charging decision epoch, of customer demand intensity, and of the battery range. The results show that the length of the charging decision epoch should allow a vehicle to be charged up to the allowed maximum energy level (80% in our case) using a fast charger. Shorter charging decision epochs would lead to higher charging frequencies and charging costs, resulting in greater customer inconvenience.
- When increasing the level of customer demand and the vehicle battery capacity, the proposed approach minimizes the charging operation time and costs to meet service needs. In contrast, the two reference charging policies apply the full-charge (80% full) policy, leading to charging more energy than necessary and resulting in higher costs and longer vehicle idle times for charging operations.
- The charging needs of other (individual-owned/commercial) EVs that compete for limited public charging facility resources are not considered here. The operator cannot know in advance the waiting time at a charger occupied by other EVs. Such a problem is well-known for EV charging at public charging stations as charging port reservation is still not available on the market [36]. In this case, vehicles can be assigned to chargers that are not occupied by other EVs in order to overcome this issue. Another alternative is to incorporate some statistical information (arrival rates and charging duration distribution) from the charging station operators to estimate the expected waiting time at a charger occupied by another non-operator-owned EV.

# Conclusions

The electrification of shared on-demand mobility services requires control over charging management as the fleet needs to frequently charge several times a day given limited public charging infrastructure. Such charging operation constraints represent significant costs for the operator due to charging queuing delays and energy costs. The operator faces the problem of scheduling the charging of the fleet in a stochastic environment with several sources of uncertainty, including the availability of charging stations, charging price variation, and stochastic customer demand. In this study, we propose a two-stage solution for handling the dynamic vehicle charging scheduling problem for dynamic dial-a-ride services using EVs that is comprised of two components: vehicle charging scheduling and vehicle–charger assignment. Charging scheduling is considered on the basis of each vehicle as a battery recharge problem, which decomposes the problem into multistage decision-making to minimize the charging costs at each stage while satisfying vehicle driving needs for the next stage. Given the charging plans of vehicles, the second component determines online vehicle–charger assignment based on the principle of the vehicle idle time for recharge minimization, considering the queuing status at the level of chargers. We apply the method to a realistic dynamic electric dial-a-ride service in Luxembourg under different charging infrastructure scenarios. The results show that significant savings can be obtained for the daily charging operations of the fleet (50 electric shuttles with 1000 customers per day): –73.4% and –76.4% in terms of the total charging waiting times, –38.1% and –39% for the total charging time, and –27.6% and –27.2% for the total charged energy costs, compared to the widely-used nearest charging station policy and minimum charging delay policy, respectively.

The approach can be extended to manage the vehicle charging scheduling of other dynamic shared mobility services such as e-taxis or ride-hailing in a dynamic environment. Future extensions can consider incorporating the charging patterns of other private EVs for more accurate waiting time estimations for chargers occupied by other private/commercial vehicles. Another direction is incorporating a day-to-day learning mechanism or a prediction model to anticipate short-term vehicle driving patterns and energy needs when demand is volatile. Moreover, the approach can be extended to consider a more realistic energy consumption model in the urban environment [37] or a mixed gasoline and electric vehicle fleet to reduce charging operational costs [26].

# Acknowledgements

The work was supported by the Luxembourg National Research Fund (C20/SC/14703944). The author thanks the Ministry of Mobility and Public Works of Luxembourg for providing the anonymized 2017 Luxmobil survey data and Francesco Viti for his helpful suggestions.

# Appendix A. Illustrative example for optimal single-vehicle battery recharge (P1)

Consider a single-vehicle battery recharge over a planning horizon from 7:00–22:00. The vehicle is fully charged at the beginning of the day. The planning horizon is divided into a sequence of charging decision epochs with a 30-minute time interval. The probability of the vehicle being in the driving state and the electricity price distribution are shown on the left side of Figure A.1. The expected waiting time to be served at charging stations over the planning horizon is depicted on the right side of Figure A.1. The battery capacity of the vehicle is assumed to be 24 kWh for the first scenario and 48 kWh for the second scenario. Two recharge policies are compared: 1) optimal single-vehicle battery recharge (P1), and 2) the on-need policy, assuming that the vehicle is recharged to 100% of its battery capacity in one epoch whenever its energy level is lower than 20%.

The parameter settings for the illustrative example are shown in Table A.1.

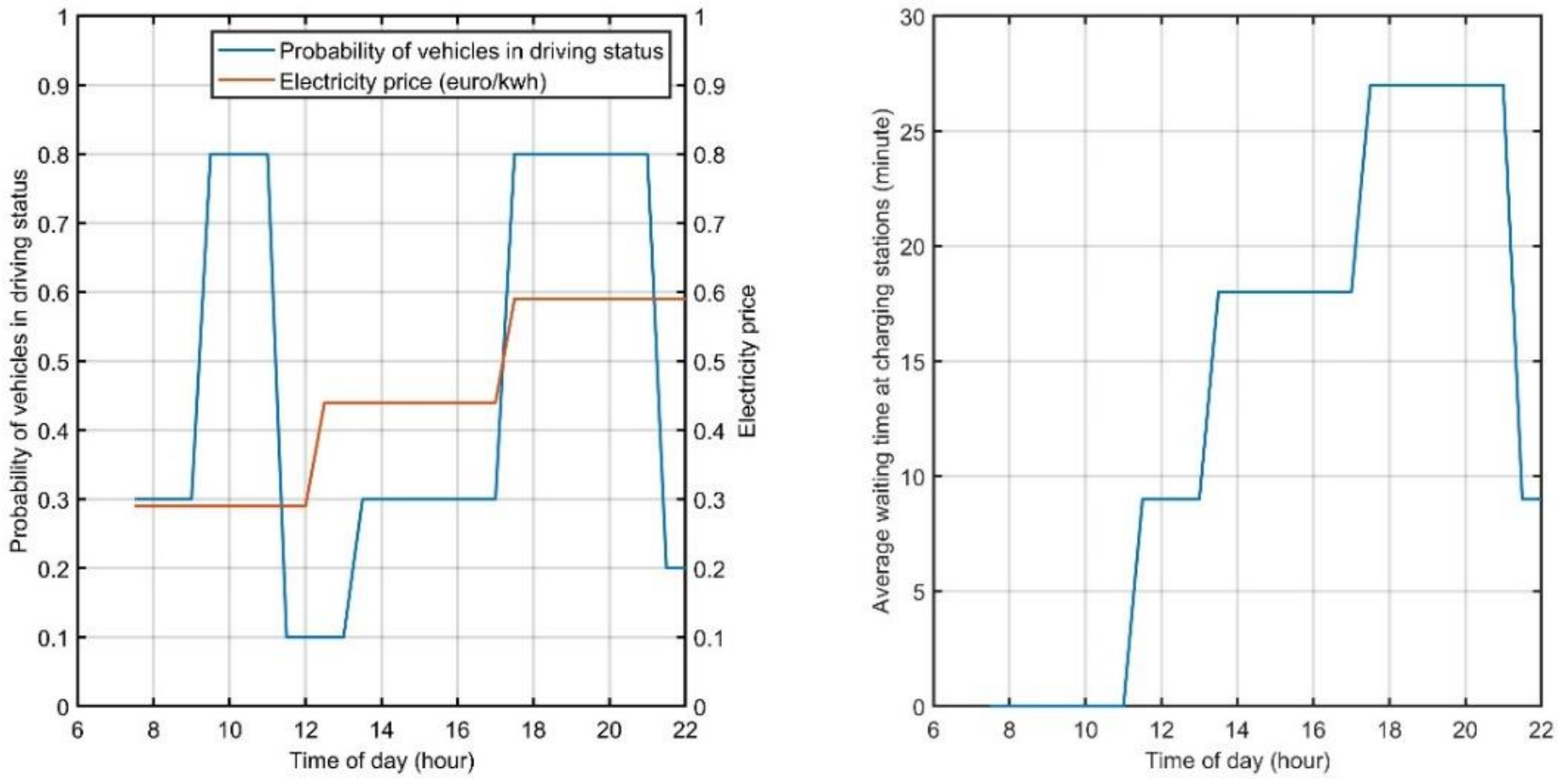

Figure A.1. Probability distributions of the vehicle being in the driving state and the electricity price (euros/kWh) over the charging decision epochs (Left); average waiting time to be served at charging stations over the charging decision epochs (Right).

Table A.1. Parameter settings for the illustrative example.

| Parameter | Value |
|---|---|
| $e_{max}$ | 24 kWh (scenario 1) and 48 kWh (scenario 2) |
| $e_{min}$ | $0.1B$ |
| $\varphi$ | 2/3 (kW/min.) |
| $\rho$ | 1/6 (euro/min.) |
| $v$ | 4/6 (km/min.) |
| $\mu$ | 0.2 (kWh/km) |
| $\bar{c}$ | 3 euros |
| $\Delta$ | 30 minutes |
| $T$ | 7:00–22:00 |

Figure A.2 (Left) shows the vehicle recharge profile under the on-need policy. For the case of 24 kWh, the vehicle needs to be recharged twice, from 14:00 to 14:30 and from 19:00 to 19:30, with around 20 kWh charged each time. For the case of 48 kWh, the vehicle goes to recharge from 19:00 to 19:30, with 42 kWh

charged. For the optimal battery recharge policy (Figure A.3), the vehicle is recharged from 11:30 to 12:00 (18 kWh charged) and from 17:30 to 18:00 (16 kWh charged) due to a lower driving probability and electricity price, given the 24 kWh battery capacity. If the battery capacity is doubled, the optimal charging time is at 11:30, with 12.7 kWh charged to meet the vehicle's driving needs until the end of the planning horizon (22:00). The battery level is always no less than $e_{min}$. Table A.2 compares the total charging costs of the two charging policies. For the optimal battery recharge policy, the total charging operation costs (Eq. (1)) are 31.28 euros (B=24 kWh) and 8.71 euros (B=48 kWh). For the on-need policy, the total costs are 40.3 euros (B=24kWh) and 36.04 euros (B=48kWh). The charging cost savings for the optimal recharge policy are 22.4% (B=24 kWh) and 75.8% (B=48 kWh). The proposed optimal charging plan allows the vehicle to be charged with minimal cost by considering the vehicle's driving needs, expected waiting times at charging stations, and opportunity cost while not available for service.

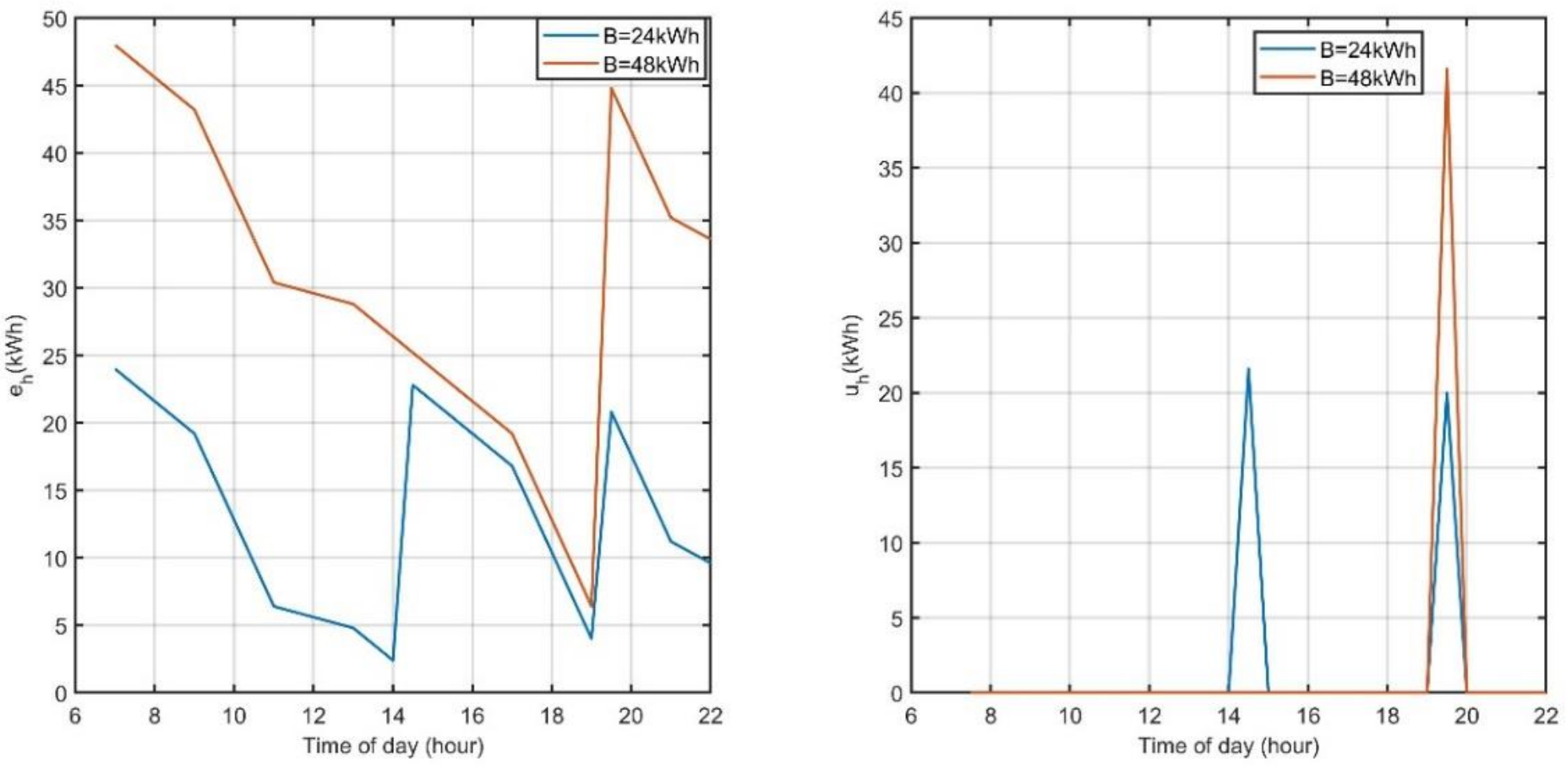

Figure A.2. Battery levels at the beginning of each epoch based on the on-need policy (Left). Recharged energy amounts in each epoch based on the on-need policy (Right).

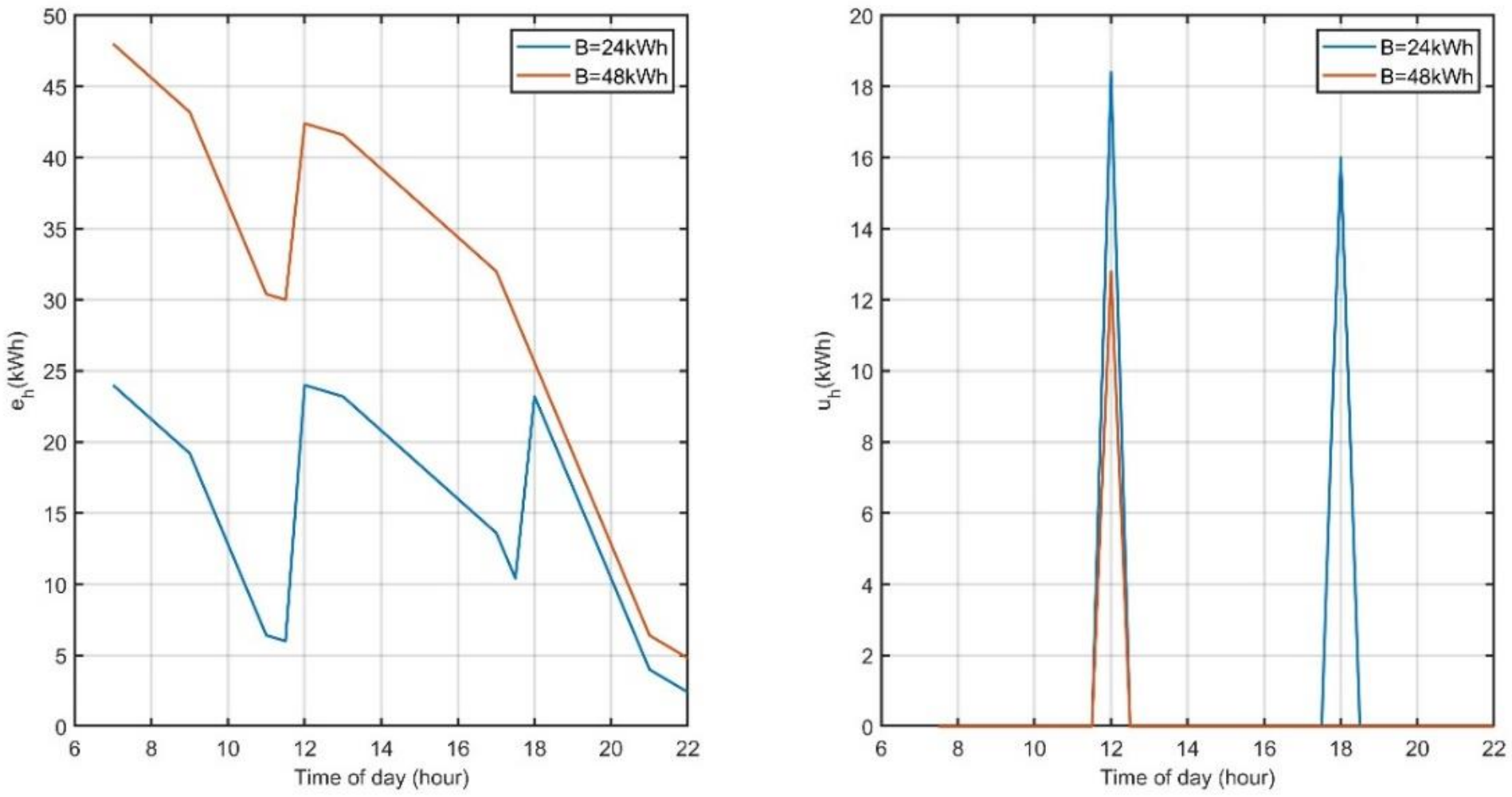

Figure A.3. Battery levels at the beginning of each epoch based on optimal battery recharge (Left). Recharged energy at each epoch based on optimal battery recharge (Right).

Table A.2. Total charging operation costs given different battery capacities.

| Battery capacity (kWh) | On-need policy (euro) | Optimal charging schedule policy (euro) | Saving |
|---|---|---|---|
| B=24 | 40.30 | 31.28 | –22.4% |
| B=48 | 36.04 | 8.71 | –75.8% |

# Appendix B. Illustrative example and computational study for optimal charging station assignment (P2)

The example is designed to illustrate the characteristics of the charging station assignment model. Consider 5 EVs with an identical battery capacity of 35.8 kWh and a full-charge driving range of 150 km, located at nodes 2, 5, 6, 7, and 9 with the same initial battery level of 20%, as shown in Figure B.1. The parameter settings for the illustrative example is shown in Table B.1. The charger power is identical at 40 kW for each charger. Vehicle speed is assumed to be 50 km/hour. The target energy levels to recharge for the vehicles, from left to right, are 80%, 40%, 50%, 80%, and 40%, respectively. There are a total of 4 chargers, of which 2 (A and B) are located at node 3 and the other two (C and D) at node 10. The available times for charger A and charger B are t=0 and t=40 min. (being occupied until t=40 min.), respectively. The available times for charger C and charger D are t=25 and t=20, respectively. The distance and travel time between any two adjacent nodes are 5 km and 6 min, respectively. As the number of vehicles is greater than that of chargers, the P2J problem is solved using the MATLAB intlinprog mixed-integer linear programming solver. The assignment results are shown in Figure B.1. The blue line represents the path that a vehicle takes to get to a charger and be charged at a higher State of Charge (SoC) level. The assignment results show that the total charging operation time is minimized with one-to-one vehicle–charger matches. Four vehicles are charged to its desired SoCs while the vehicle on the left side is not charged due to the charger availability limit (4 chargers vs. 5 vehicles) and its higher target charging level (80%). This demonstrates that when the number of vehicles to be charged is greater than that of chargers, the vehicles with a lower energy demand are prioritized for charging when all things are equal. This allows more vehicles to be available earlier to serve customers. Table B.2 shows the detailed result of the assignment model with the obtained value of the objective function Z*= 182.92. The arrival time, waiting time, and charging time of each vehicle are reported in Table B.2.

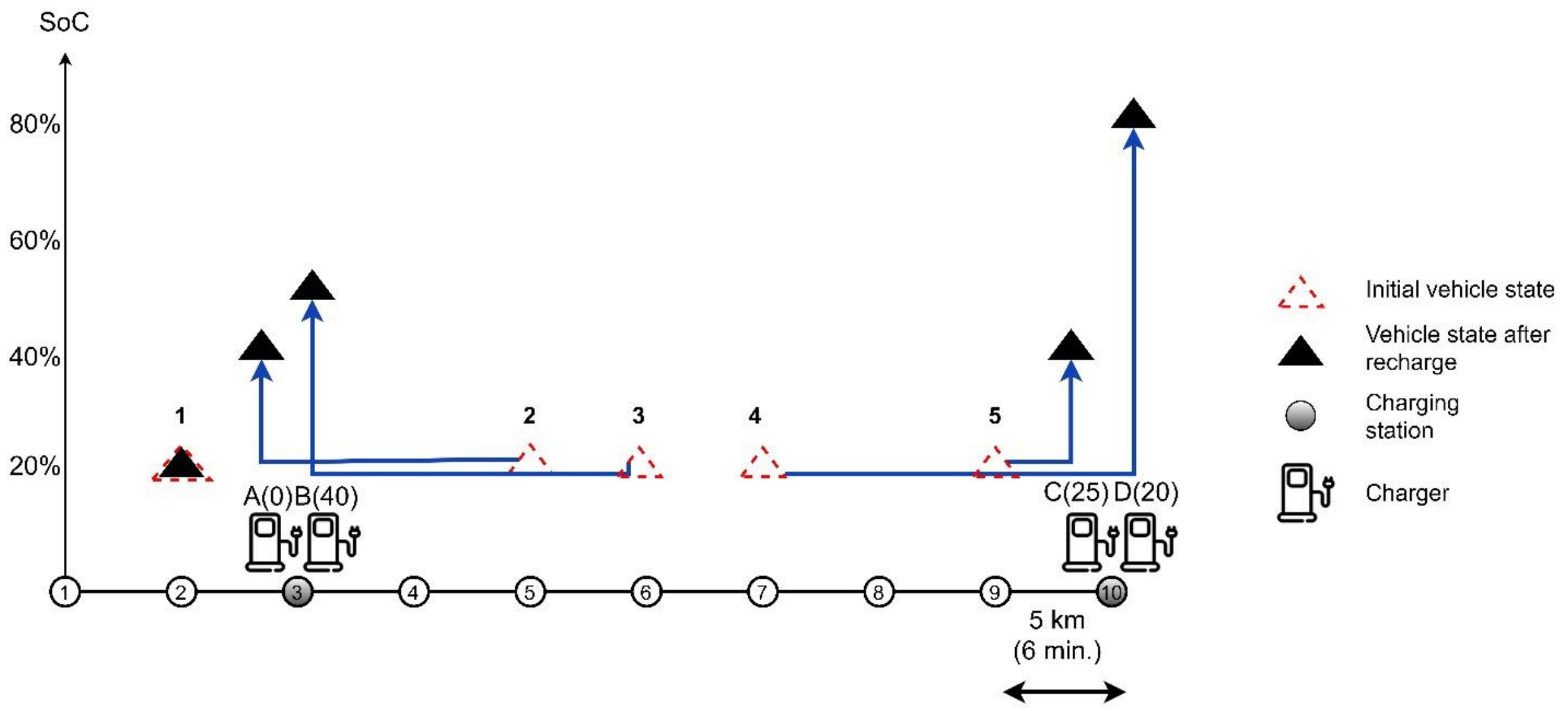

Figure B.1. Illustrative example of the charging station assignment model with 5 vehicles and 4 chargers (A,B,C, and D).

Table B.1. Parameter settings for the illustrative example.

| | |
|---|---|
| $B$ (battery capacity) | 35.8 (kWh) |
| $e_{min}$ | $0.1B$ |
| Driving range (full charge) | 150 (km) |
| $\theta_1, \theta_2$ | 1 |
| $\varphi$ (charging rate) | 40 (kW/hour) |
| Location of chargers | (3,3,10,10) |
| Location of vehicles | (2,5,6,7,9) |
| $e_i$ | (0.2,0.2,0.2,0.2,0.2) *B (kWh) |
| Vehicle speed | 5/6 (km/min.) |
| $\varphi_j, j = 1, \dots, 4$ | 0.2386667 (kWh/km) |
| $t_j^A, j = 1, \dots, 4$ | (0,40,25,20) (min.) |
| $e_i^*, i = 1, \dots, 5$ | (0.8,0.4,0.5,0.8,0.4)*B |

Table B.2. Results of optimal charging assignment for the illustrative example.

| Vehicle | SoC (%of B) after recharge | Assigned charger | Arrival time | Waiting time | Charging time |
|---|---|---|---|---|---|
| 1 | 20 | - | - | - | - |
| 2 | 40 | A | 12.0 | 0.0 | 14.3 |
| 3 | 50 | B | 18.0 | 22.0 | 21.5 |
| 4 | 80 | D | 18.0 | 2.0 | 37.6 |
| 5 | 40 | C | 6.0 | 19.0 | 12.5 |
| Z* | 182.92 | | | | |

Remark: B is the battery capacity. Time is measured in minutes.

For larger problems, we generate two sets of test instances; each has 9 subsets of problems with different sizes of $|I|$ and $|J|$, and each subset has 3 randomly generated test instances. The first set of test instances is related to the P2 problem ($|I| \leq |J|$). The second set is related to the P2J problem ($|I| > |J|$). Note that for the second set of instances, the number of chargers $|J|$ is randomly generated given $|J| < |I|$. The locations of vehicles and chargers are randomly generated within a rectangular area within ($[-50,50] \times [-50,50]$). A vehicle's initial energy level $e_i$, target energy level $e_i^*$ after recharge, and a chargers' available time $t_j^A$ are randomly generated based on the parameters shown in Table B.3. The performance of the proposed heuristic is compared with the exact solution obtained by the MATLAB mixed-integer linear programming solver. In terms of the P2 problem size, the number of constraints in terms of Eqs. (11)–(16) is $|I| + |J| + 4|I||J|$. The test instances are publicly available at https://github.com/tym2021.

Table B.3. Parameters for test instances generation for the Lagrangian relaxation algorithm.

| Variable | Value |
|---|---|
| $B$ | 35.8 kWh |
| $e_{min}$ | $0.1B$ |
| $\|I\|$ | 10,20,30,40,50, 100,200,400, and 1000 |
| $\|J\|$ | 10,20,30,40,50, 100,200,400, and 1000 |
| Location of chargers/vehciles | rand($[-50,50] \times [-50,50]$) |
| $e_i$ | rand(0.4,0.5)*B (kWh) |
| $e_i^*$ | rand (0.7,1) *B (kWh) |
| $t_j^A$ | rand (0,30) (min.) |
| $\theta_1, \theta_2$ | 1 |
| $\varphi$ (charging rate) | 40 (kWh) |
| $\varphi_j, j = 1, \dots, 4$ | 0.2387 (kWh/km) |

The maximum number of iterations of the LR algorithm is set as 2000. Different values of the step size adjustment constant $0 < \delta \leq 2$ are tested and, finally, we pick 0.6 because it shows the best performance in most cases. The gap tolerance $\bar{\varepsilon}$ is set as 0.0001 for |I|<1000 and 0.005 otherwise. The computational results of the LR algorithm are shown in Table B.4. The results for the P2 problems are reported in the left column, while the right column reports the results for the P2J problems. The LR algorithm finds near-optimal solutions with an optimality gap of 0.5% in less than 3 minutes for the test instances of 1000 vehicles in the P2 and P2J problems. For the test instances of 200 vehicles, the optimality gap is around 0.16% for the P2 problem and 0.09% for the P2J problem, both found in around 30 seconds. However, the commercial exact solution solver cannot find solutions within 1 hour for |I|=50 and more. Our computational study shows that the proposed LR algorithm is suitable for large-scale real-time application for EV charging station assignment.

Table B.4. Computational results of the Lagrangian relaxation algorithm for the P2 and P2J problems.

| Prob. | P2 | | | | | P2J | | | | |
|---|---|---|---|---|---|---|---|---|---|---|
| | $\|I\|$ | $\|J\|$ | Gap | CPU time (seconds) | | $\|I\|$ | $\|J\|$ (avg.) | Gap | CPU time (seconds) | |
| | | | | LR | Exact | | | | LR | Exact |
| 1 | 10 | 10 | 0 | 0.1 | 0.2 | 10 | 4.7 | 0.00% | 0.1 | 0.4 |
| 2 | 20 | 20 | 0.03% | 2.9 | 0.4 | 20 | 5.0 | 0.01% | 0.3 | 0.3 |
| 3 | 30 | 30 | 0.04% | 3.9 | 0.4 | 30 | 14.3 | 0.01% | 1.7 | 0.3 |
| 4 | 40 | 40 | 0.09% | 5.3 | 130.8 | 40 | 15.0 | 0.02% | 2.3 | 1.7 |
| 5 | 50 | 50 | 0.11% | 5.7 | NA | 50 | 37.3 | 0.03% | 4.5 | NA |
| 6 | 100 | 100 | 0.09% | 11.8 | NA | 100 | 51.3 | 0.03% | 12.5 | NA |
| 7 | 200 | 200 | 0.16% | 25.5 | NA | 200 | 117.0 | 0.09% | 36.5 | NA |
| 8 | 400 | 400 | 0.26% | 121.8 | NA | 400 | 104.3 | 0.06% | 115.4 | NA |
| 9 | 1000 | 1000 | 0.50% | 168.3 | NA | 1000 | 522.3 | 0.50% | 170.5 | NA |

Remark: 1. GAP $=\frac{z_{UB}-z_{LB}}{z_{UB}}$. 2. NA means the exact solution cannot be found given a one-hour computation time. 3. The reported results are based on the average of three randomly generated test instances for each problem size.